\documentclass[12pt]{article}
\usepackage{amsmath,latexsym,amssymb}
\begin{document}
\headsep 0.5 true cm

\begin{center}
{\Large\bf $E_8$-singularity, invariant theory and modular forms}
\vskip 1.0 cm
{\large\bf Lei Yang}
\end{center}
\vskip 1.5 cm

\begin{center}
{\large\bf Abstract}
\end{center}

\vskip 0.5 cm

  As an algebraic surface, the equation of $E_8$-singularity
$x^5+y^3+z^2=0$ can be obtained from a quotient $C_Y/\text{SL}(2, 13)$
over the modular curve $X(13)$, where $Y \subset \mathbb{CP}^5$ is a
complete intersection curve given by a system of $\text{SL}(2, 13)$-invariant
polynomials and $C_Y$ is a cone over $Y$. It is different from the
Kleinian singularity $\mathbb{C}^2/\Gamma$, where $\Gamma$ is the
binary icosahedral group. This gives a negative answer to Arnol'd
and Brieskorn's questions about the mysterious relation between the
icosahedron and $E_8$, i.e., the $E_8$-singularity is not necessarily
the Kleinian icosahedral singularity. In particular, the equation of
$E_8$-singularity possesses infinitely many kinds of distinct modular
parametrizations, and there are infinitely many kinds of distinct
constructions of the $E_8$-singularity. They form a variation of the
$E_8$-singularity structure over the modular curve $X(13)$, for which
we give its algebraic version, geometric version, $j$-function version
and the version of Poincar\'{e} homology $3$-sphere as well as its higher
dimensional lifting, i.e., Milnor's exotic $7$-sphere. Moreover, there
are variations of $Q_{18}$ and $E_{20}$-singularity structures over $X(13)$.
Thus, three different algebraic surfaces, the equations of $E_8$, $Q_{18}$
and $E_{20}$-singularities can be realized from the same quotients
$C_Y/\text{SL}(2, 13)$ over the modular curve $X(13)$ and have the same
modular parametrizations.

\vskip 5.0 cm

\begin{center}
{\large\bf Contents}
\end{center}

$$\aligned
  &1. \text{\quad Introduction}\\
  &2. \text{\quad Standard construction: $E_8$-singularity
                  as the icosahedral singularity}\\
  &3. \text{\quad Invariant theory and modular forms for
                  $\text{SL}(2, 13)$}\\
  &4. \text{\quad A different construction: $E_8$-singularity
                  from $C_Y/\text{SL}(2, 13)$ over $X(13)$}\\
  &\quad \quad \text{and a variation of the $E_8$-singularity
                  structure over $X(13)$}\\
  &5. \text{\quad $Q_{18}$ and $E_{20}$-singularities from
                  $C_Y/\text{SL}(2, 13)$ over $X(13)$ and}\\
  &\quad \quad \text{variations of $Q_{18}$ and $E_{20}$-singularity
                  structures over $X(13)$}
\endaligned$$

\vskip 2.0 cm

\begin{center}
{\large\bf 1. Introduction}
\end{center}

  By a rational double point or a simple singularity we understand the
singularity of the quotient of $\mathbb{C}^2$ by the action of a finite
subgroup of $\text{SL}(2, \mathbb{C})$ (see \cite{Sl1}). Let $\Gamma$ be
a finite subgroup of $\text{SL}(2, \mathbb{C})$. Then $\Gamma$ is one of
the following: a cyclic group of order $\ell \geq 1$ $(A_{\ell})$, a
binary dihedral group of order $4(\ell-2)$, $\ell \geq 4$ $(D_{\ell})$,
the binary tetrahedral group $(E_6)$, the binary octahedral group $(E_7)$,
or the binary icosahedral group $(E_8)$. In 1874, Klein showed that the
ring of polynomials in two variables which are invariant under $\Gamma$
is generated by three elements $x$, $y$ and $z$, which satisfy the
following relation
$$\aligned
 &A_{\ell \geq 1} \quad  & x^{\ell+1}+y^2+z^2=0,\\
 &D_{\ell \geq 4} \quad  & x^{\ell-1}+x y^2+z^2=0,\\
 &E_6             \quad  & x^4+y^3+z^2=0,\\
 &E_7             \quad  & x^3 y+y^3+z^2=0,\\
 &E_8             \quad  & x^5+y^3+z^2=0.
\endaligned$$
These results of Klein on the invariant theory of the binary polyhedral
groups were a starting point for later developments. In the minimal
resolution of such a singularity an intersection configuration of the
components of the exceptional divisor appears which can be described in
a simple way by a Dynkin diagram of type $A_{\ell}$, $D_{\ell}$, $E_6$,
$E_7$ or $E_8$. Up to analytic isomorphism, these diagrams classify the
corresponding singularities. In other words, the $ADE$ singularities are
the Kleinian singularities, i.e., the quotient singularities of $\mathbb{C}^2$
by a finite subgroup of $\text{SL}(2, \mathbb{C})$. In particular, the
$E_8$-singularity is the icosahedral singularity $\mathbb{C}^2/\Gamma$,
where $\Gamma$ is the binary icosahedral group.

  In the present paper, we will show that the $E_8$-singularity can be
obtained from a quotient $C_Y/\text{SL}(2, 13)$ over the modular curve
$X(13)$, where $Y \subset \mathbb{CP}^5$ is a complete intersection curve
given by a system of $\text{SL}(2, 13)$-invariant polynomials and $C_Y$
is a cone over $Y$. This gives an infinitely many kinds of distinct
constructions of the $E_8$-singularity which are different from the
icosahedral singularity. They form a variation of the $E_8$-singularity
structure over the modular curve $X(13)$, for which we give its algebraic
version, geometric version, $j$-function version and the version of
Poincar\'{e} homology $3$-sphere as well as its higher dimensional lifting,
i.e., Milnor's exotic $7$-sphere. Our construction is based on the invariant
theory for the group $\text{SL}(2, 13)$. We obtain some invariants
$\Phi_{12}$, $\Phi_{20}$ and $\Phi_{30}$ for $\text{SL}(2, 13)$. Furthermore,
over the modular curve $X(13)$, these invariants are modular forms which
satisfy the equation of the $E_8$-singularity. Hence, we obtain an
homomorphism from the ring
$$\mathbb{C}[\Phi_{12}, \Phi_{20}, \Phi_{30}]/
  (\Phi_{20}^3-\Phi_{30}^2-1728 \Phi_{12}^5)$$
to the ring of invariants
$$\left[\mathbb{C}[z_1, z_2, z_3, z_4, z_5, z_6]/I \right]^{\text{SL}(2, 13)}$$
of the group $\text{SL}(2, 13)$ acting on the cone $C_Y \subset \mathbb{C}^6$
over the modular curve $X(13)$, where $I$ is an ideal generated by a system
of algebraic relations satisfied by the variables $z_1$, $\ldots$ $z_6$, $Y$
is the algebraic curve corresponding the ideal $I$.

  Let us begin with the invariant theory for $\text{SL}(2, 13)$. The
representation of $\text{SL}(2, 13)$ we will consider is the unique
six-dimensional irreducible complex representation for which the
eigenvalues of $\left(\begin{matrix} 1 & 1\\ 0 & 1 \end{matrix}\right)$
are the $\exp (\underline{a} . 2 \pi i/13)$ for $\underline{a}$ a
non-square mod $13$. We will give an explicit realization of this
representation. This explicit realization will play a major role for
giving a complete system of invariants associated to $\text{SL}(2, 13)$.
Recall that the six-dimensional representation of the finite group
$\text{SL}(2, 13)$ of order $2184$, which acts on the five-dimensional
projective space
$$\mathbb{CP}^5=\{ (z_1, z_2, z_3, z_4, z_5, z_6): z_i \in \mathbb{C}
  \quad (i=1, 2, 3, 4, 5, 6) \}.$$
This representation is defined over the cyclotomic field
$\mathbb{Q}(e^{\frac{2 \pi i}{13}})$. Put
$$S=-\frac{1}{\sqrt{13}} \begin{pmatrix}
  \zeta^{12}-\zeta & \zeta^{10}-\zeta^3 & \zeta^4-\zeta^9
& \zeta^5-\zeta^8 & \zeta^2-\zeta^{11} & \zeta^6-\zeta^7\\
  \zeta^{10}-\zeta^3 & \zeta^4-\zeta^9 & \zeta^{12}-\zeta
& \zeta^2-\zeta^{11} & \zeta^6-\zeta^7 & \zeta^5-\zeta^8\\
  \zeta^4-\zeta^9 & \zeta^{12}-\zeta & \zeta^{10}-\zeta^3
& \zeta^6-\zeta^7 & \zeta^5-\zeta^8 & \zeta^2-\zeta^{11}\\
  \zeta^5-\zeta^8 & \zeta^2-\zeta^{11} & \zeta^6-\zeta^7
& \zeta-\zeta^{12} & \zeta^3-\zeta^{10} & \zeta^9-\zeta^4\\
  \zeta^2-\zeta^{11} & \zeta^6-\zeta^7 & \zeta^5-\zeta^8
& \zeta^3-\zeta^{10} & \zeta^9-\zeta^4 & \zeta-\zeta^{12}\\
  \zeta^6-\zeta^7 & \zeta^5-\zeta^8 & \zeta^2-\zeta^{11}
& \zeta^9-\zeta^4 & \zeta-\zeta^{12} & \zeta^3-\zeta^{10}
\end{pmatrix}$$
and
$$T=\text{diag}(\zeta^7, \zeta^{11}, \zeta^8, \zeta^6, \zeta^2, \zeta^5)$$
where $\zeta=\exp(2 \pi i/13)$. We have
$$S^2=-I, \quad T^{13}=(ST)^3=I.\eqno{(1.1)}$$
Let
$G=\langle S, T \rangle$, then $G \cong \text{SL}(2, 13)$.
We construct some $G$-invariant polynomials in six variables
$z_1, \ldots, z_6$. Let
$$w_{\infty}=13 \mathbf{A}_0^2, \quad
  w_{\nu}=(\mathbf{A}_0+\zeta^{\nu} \mathbf{A}_1+\zeta^{4 \nu} \mathbf{A}_2
  +\zeta^{9 \nu} \mathbf{A}_3+\zeta^{3 \nu} \mathbf{A}_4+\zeta^{12 \nu}
  \mathbf{A}_5+\zeta^{10 \nu} \mathbf{A}_6)^2\eqno{(1.2)}$$
for $\nu=0, 1, \ldots, 12$, where the senary quadratic forms (quadratic
forms in six variables) $\mathbf{A}_j$ $(j=0, 1, \ldots, 6)$ are given by
$$\left\{\aligned
  \mathbf{A}_0 &=z_1 z_4+z_2 z_5+z_3 z_6,\\
  \mathbf{A}_1 &=z_1^2-2 z_3 z_4,\\
  \mathbf{A}_2 &=-z_5^2-2 z_2 z_4,\\
  \mathbf{A}_3 &=z_2^2-2 z_1 z_5,\\
  \mathbf{A}_4 &=z_3^2-2 z_2 z_6,\\
  \mathbf{A}_5 &=-z_4^2-2 z_1 z_6,\\
  \mathbf{A}_6 &=-z_6^2-2 z_3 z_5.
\endaligned\right.\eqno{(1.3)}$$
Then $w_{\infty}$, $w_{\nu}$ for $\nu=0, \ldots, 12$ are the roots of a
polynomial of degree fourteen. The corresponding equation is just the
Jacobian equation of degree fourteen (see \cite{K}, pp.161-162). On the
other hand, set
$$\delta_{\infty}=13^2 \mathbf{G}_0, \quad
  \delta_{\nu}=-13 \mathbf{G}_0+\zeta^{\nu} \mathbf{G}_1+\zeta^{2 \nu}
  \mathbf{G}_2+\cdots+\zeta^{12 \nu} \mathbf{G}_{12}\eqno{(1.4)}$$
for $\nu=0, 1, \ldots, 12$, where the senary sextic forms (i.e., sextic
forms in six variables) $\mathbf{G}_j$ $(j=0, 1, \ldots, 12)$ are given
by
$$\left\{\aligned
  \mathbf{G}_0 =&\mathbf{D}_0^2+\mathbf{D}_{\infty}^2,\\
  \mathbf{G}_1 =&-\mathbf{D}_7^2+2 \mathbf{D}_0 \mathbf{D}_1+10 \mathbf{D}_{\infty}
                 \mathbf{D}_1+2 \mathbf{D}_2 \mathbf{D}_{12}+\\
                &-2 \mathbf{D}_3 \mathbf{D}_{11}-4 \mathbf{D}_4 \mathbf{D}_{10}
                 -2 \mathbf{D}_9 \mathbf{D}_5,\\
  \mathbf{G}_2 =&-2 \mathbf{D}_1^2-4 \mathbf{D}_0 \mathbf{D}_2+6 \mathbf{D}_{\infty}
                 \mathbf{D}_2-2 \mathbf{D}_4 \mathbf{D}_{11}+\\
                &+2 \mathbf{D}_5 \mathbf{D}_{10}-2 \mathbf{D}_6 \mathbf{D}_9-2
                 \mathbf{D}_7 \mathbf{D}_8,\\
  \mathbf{G}_3 =&-\mathbf{D}_8^2+2 \mathbf{D}_0 \mathbf{D}_3+10 \mathbf{D}_{\infty}
                 \mathbf{D}_3+2 \mathbf{D}_6 \mathbf{D}_{10}+\\
                &-2 \mathbf{D}_9 \mathbf{D}_7-4 \mathbf{D}_{12} \mathbf{D}_4
                 -2 \mathbf{D}_1 \mathbf{D}_2,\\
  \mathbf{G}_4 =&-\mathbf{D}_2^2+10 \mathbf{D}_0 \mathbf{D}_4-2 \mathbf{D}_{\infty}
                 \mathbf{D}_4+2 \mathbf{D}_5 \mathbf{D}_{12}+\\
                &-2 \mathbf{D}_9 \mathbf{D}_8-4 \mathbf{D}_1 \mathbf{D}_3-2
                 \mathbf{D}_{10} \mathbf{D}_7,\\
  \mathbf{G}_5 =&-2 \mathbf{D}_9^2-4 \mathbf{D}_0 \mathbf{D}_5+6 \mathbf{D}_{\infty}
                 \mathbf{D}_5-2 \mathbf{D}_{10} \mathbf{D}_8+\\
                &+2 \mathbf{D}_6 \mathbf{D}_{12}-2 \mathbf{D}_2 \mathbf{D}_3
                 -2 \mathbf{D}_{11} \mathbf{D}_7,\\
  \mathbf{G}_6 =&-2 \mathbf{D}_3^2-4 \mathbf{D}_0 \mathbf{D}_6+6 \mathbf{D}_{\infty}
                 \mathbf{D}_6-2 \mathbf{D}_{12} \mathbf{D}_7+\\
                &+2 \mathbf{D}_2 \mathbf{D}_4-2 \mathbf{D}_5 \mathbf{D}_1-2
                 \mathbf{D}_8 \mathbf{D}_{11},\\
  \mathbf{G}_7 =&-2 \mathbf{D}_{10}^2+6 \mathbf{D}_0 \mathbf{D}_7+4 \mathbf{D}_{\infty}
                 \mathbf{D}_7-2 \mathbf{D}_1 \mathbf{D}_6+\\
                &-2 \mathbf{D}_2 \mathbf{D}_5-2 \mathbf{D}_8 \mathbf{D}_{12}-2
                 \mathbf{D}_9 \mathbf{D}_{11},\\
  \mathbf{G}_8 =&-2 \mathbf{D}_4^2+6 \mathbf{D}_0 \mathbf{D}_8+4 \mathbf{D}_{\infty}
                 \mathbf{D}_8-2 \mathbf{D}_3 \mathbf{D}_5+\\
                &-2 \mathbf{D}_6 \mathbf{D}_2-2 \mathbf{D}_{11} \mathbf{D}_{10}-2
                 \mathbf{D}_1 \mathbf{D}_7,\\
  \mathbf{G}_9 =&-\mathbf{D}_{11}^2+2 \mathbf{D}_0 \mathbf{D}_9+10 \mathbf{D}_{\infty}
                 \mathbf{D}_9+2 \mathbf{D}_5 \mathbf{D}_4+\\
                &-2 \mathbf{D}_1 \mathbf{D}_8-4 \mathbf{D}_{10} \mathbf{D}_{12}-2
                 \mathbf{D}_3 \mathbf{D}_6,\\
  \mathbf{G}_{10} =&-\mathbf{D}_5^2+10 \mathbf{D}_0 \mathbf{D}_{10}-2 \mathbf{D}_{\infty}
                    \mathbf{D}_{10}+2 \mathbf{D}_6 \mathbf{D}_4+\\
                   &-2 \mathbf{D}_3 \mathbf{D}_7-4 \mathbf{D}_9 \mathbf{D}_1-2
                    \mathbf{D}_{12} \mathbf{D}_{11},\\
  \mathbf{G}_{11} =&-2 \mathbf{D}_{12}^2+6 \mathbf{D}_0 \mathbf{D}_{11}+4 \mathbf{D}_{\infty}
                    \mathbf{D}_{11}-2 \mathbf{D}_9 \mathbf{D}_2+\\
                   &-2 \mathbf{D}_5 \mathbf{D}_6-2 \mathbf{D}_7 \mathbf{D}_4-2
                    \mathbf{D}_3 \mathbf{D}_8,\\
  \mathbf{G}_{12} =&-\mathbf{D}_6^2+10 \mathbf{D}_0 \mathbf{D}_{12}-2 \mathbf{D}_{\infty}
                    \mathbf{D}_{12}+2 \mathbf{D}_2 \mathbf{D}_{10}+\\
                   &-2 \mathbf{D}_1 \mathbf{D}_{11}-4 \mathbf{D}_3 \mathbf{D}_9-2
                    \mathbf{D}_4 \mathbf{D}_8.
\endaligned\right.\eqno{(1.5)}$$
Here, the senary cubic forms (cubic forms in six variables)
$\mathbf{D}_j$ $(j=0$, $1$, $\ldots$, $12$, $\infty)$ are given as follows:
$$\left\{\aligned
  \mathbf{D}_0 &=z_1 z_2 z_3,\\
  \mathbf{D}_1 &=2 z_2 z_3^2+z_2^2 z_6-z_4^2 z_5+z_1 z_5 z_6,\\
  \mathbf{D}_2 &=-z_6^3+z_2^2 z_4-2 z_2 z_5^2+z_1 z_4 z_5+3 z_3 z_5 z_6,\\
  \mathbf{D}_3 &=2 z_1 z_2^2+z_1^2 z_5-z_4 z_6^2+z_3 z_4 z_5,\\
  \mathbf{D}_4 &=-z_2^2 z_3+z_1 z_6^2-2 z_4^2 z_6-z_1 z_3 z_5,\\
  \mathbf{D}_5 &=-z_4^3+z_3^2 z_5-2 z_3 z_6^2+z_2 z_5 z_6+3 z_1 z_4 z_6,\\
  \mathbf{D}_6 &=-z_5^3+z_1^2 z_6-2 z_1 z_4^2+z_3 z_4 z_6+3 z_2 z_4 z_5,\\
  \mathbf{D}_7 &=-z_2^3+z_3 z_4^2-z_1 z_3 z_6-3 z_1 z_2 z_5+2 z_1^2 z_4,\\
  \mathbf{D}_8 &=-z_1^3+z_2 z_6^2-z_2 z_3 z_5-3 z_1 z_3 z_4+2 z_3^2 z_6,\\
  \mathbf{D}_9 &=2 z_1^2 z_3+z_3^2 z_4-z_5^2 z_6+z_2 z_4 z_6,\\
  \mathbf{D}_{10} &=-z_1 z_3^2+z_2 z_4^2-2 z_4 z_5^2-z_1 z_2 z_6,\\
  \mathbf{D}_{11} &=-z_3^3+z_1 z_5^2-z_1 z_2 z_4-3 z_2 z_3 z_6+2 z_2^2 z_5,\\
  \mathbf{D}_{12} &=-z_1^2 z_2+z_3 z_5^2-2 z_5 z_6^2-z_2 z_3 z_4,\\
  \mathbf{D}_{\infty}&=z_4 z_5 z_6.
\endaligned\right.\eqno{(1.6)}$$
Then $\delta_{\infty}$, $\delta_{\nu}$ for $\nu=0, \ldots, 12$ are the
roots of a polynomial of degree fourteen. The corresponding equation is
not the Jacobian equation. Let $S_d=\mathbb{C}[\Phi_{m, n}]$ be the
subalgebra of $\mathbb{C}[z_1, \ldots, z_6]$ generated by the invariant
homogeneous polynomials $\Phi_{m, n}$ given by
$$\Phi_{m, n}=w_0^m \delta_0^n+w_1^m \delta_1^n+\cdots+w_{12}^m \delta_{12}^n+
              w_{\infty}^m \delta_{\infty}^n,\eqno{(1.7)}$$
with degree $d=4m+6n$. Let
$$S=\bigoplus_{d \geq 0} S_d.$$
Then
$$S \subseteq \mathbb{C}[\mathbf{A}_0, \mathbf{A}_1, \ldots, \mathbf{A}_6,
  \mathbf{G}_0, \mathbf{G}_1, \ldots, \mathbf{G}_{12}]^G \subseteq
  \mathbb{C}[z_1, \ldots, z_6]^G.$$

  Let $x_i(z)=\eta(z) a_i(z)$ $(1 \leq i \leq 6)$, where
$$\left\{\aligned
  a_1(z) &:=e^{-\frac{11 \pi i}{26}} \theta \begin{bmatrix}
            \frac{11}{13}\\ 1 \end{bmatrix}(0, 13z),\\
  a_2(z) &:=e^{-\frac{7 \pi i}{26}} \theta \begin{bmatrix}
            \frac{7}{13}\\ 1 \end{bmatrix}(0, 13z),\\
  a_3(z) &:=e^{-\frac{5 \pi i}{26}} \theta \begin{bmatrix}
            \frac{5}{13}\\ 1 \end{bmatrix}(0, 13z),\\
  a_4(z) &:=-e^{-\frac{3 \pi i}{26}} \theta \begin{bmatrix}
            \frac{3}{13}\\ 1 \end{bmatrix}(0, 13z),\\
  a_5(z) &:=e^{-\frac{9 \pi i}{26}} \theta \begin{bmatrix}
            \frac{9}{13}\\ 1 \end{bmatrix}(0, 13z),\\
  a_6(z) &:=e^{-\frac{\pi i}{26}} \theta \begin{bmatrix}
            \frac{1}{13}\\ 1 \end{bmatrix}(0, 13z)
\endaligned\right.\eqno{(1.8)}$$
are theta constants of order $13$ and
$\eta(z):=q^{\frac{1}{24}} \prod_{n=1}^{\infty} (1-q^n)$ with
$q=e^{2 \pi i z}$ is the Dedekind eta function which are all
defined in the upper-half plane
$\mathbb{H}=\{ z \in \mathbb{C}: \text{Im}(z)>0 \}$.
In fact, the weight of $x_i(z)$ is $1$ and the parabolic modular
forms $a_i(z)$ of weight $\frac{1}{2}$ given by (1.8) form a
multiplier-system in the sense of the following (see (3.16) in
Proposition 3.2):
$$\mathbf{A}(z+1)=e^{-\frac{3 \pi i}{4}} T \mathbf{A}(z), \quad
  \mathbf{A}\left(-\frac{1}{z}\right)=e^{\frac{\pi i}{4}} \sqrt{z}
  S \mathbf{A}(z),\eqno{(1.9)}$$
where $S$ and $T$ are given as above, $0<\text{arg} \sqrt{z} \leq \pi/2$
and
$$\mathbf{A}(z):=(a_1(z), a_2(z), a_3(z), a_4(z), a_5(z), a_6(z))^{T}.\eqno{(1.10)}$$
We will show that there is a morphism
$$\Phi: X \to Y \subset \mathbb{CP}^5\eqno{(1.11)}$$
with $\Phi(z)=(x_1(z), \ldots, x_6(z))$, where $X=X(13)$ is the modular
curve $\overline{\Gamma(13) \backslash \mathbb{H}}$ and $Y$ is a complete
intersection algebraic curve with multi-degree $(4, 8, 10, 14)$ corresponding
to the ideal
$$I=I(Y)=(\Phi_4, \Phi_8, \Phi_{10}, \Phi_{14}),\eqno{(1.12)}$$
where
$$\Phi_4=\Phi_{1, 0}, \quad \Phi_8=\Phi_{2, 0}, \quad
  \Phi_{10}=\Phi_{1, 1}, \quad \Phi_{14}=\Phi_{2, 1}.\eqno{(1.13)}$$
Each $\Phi_i$ $(i=4, 8, 10, 14)$ corresponds to a unique $\Phi_{m, n}$ with
degree $i=4m+6n$. The significance of the algebraic curve $Y$ is that the
finite group $G$ acts linearly on $\mathbb{C}^6$ and on $\mathbb{CP}^5$
leaving invariant $Y \subset \mathbb{CP}^5$ and the cone
$C_Y \subset \mathbb{C}^6$. Moreover, it is (1.9) that gives an explicit
realization of the isomorphism between the unique sub-representation of
parabolic modular forms of weight $\frac{1}{2}$ on $X(13)$ and the above
six-dimensional complex representation of $\text{SL}(2, 13)$ generated by
$S$ and $T$. Our main theorems are the following:

\textbf{Theorem 1.1.} (Main Theorem 1) (A variation of the $E_8$-singularity
structure over the modular curve $X$: algebraic version) {\it The equation of
$E_8$-singularity $$\Phi_{20}^3-\Phi_{30}^2-1728 \Phi_{12}^5=0$$
possesses an infinitely many kinds of distinct modular parametrizations
$($with the cardinality of the continuum in ZFC set theory$)$
$$(\Phi_{12}, \Phi_{20}, \Phi_{30})=(\Phi_{12}^{\lambda},
   \Phi_{20}^{\mu}, \Phi_{30}^{\gamma})\eqno{(1.14)}$$
over the modular curve $X$ as follows$:$
$$\left\{\aligned
  \Phi_{12}^{\lambda} &=\lambda \Phi_{3, 0}+(1-\lambda) \Phi_{0, 2}
              \quad \text{mod $\mathfrak{a}_1$},\\
  \Phi_{20}^{\mu} &=\mu \Phi_{5, 0}+(1-\mu) \Phi_{2, 2}
              \quad \text{mod $\mathfrak{a}_2$},\\
  \Phi_{30}^{\gamma} &=\gamma_1 \Phi_{0, 5}+\gamma_2 \Phi_{3, 3}+
              (1-\gamma_1-\gamma_2) \Phi_{6, 1} \quad \text{mod $\mathfrak{a}_3$},
\endaligned\right.\eqno{(1.15)}$$
where $\Phi_{12}$, $\Phi_{20}$ and $\Phi_{30}$ are invariant homogeneous
polynomials of degree $12$, $20$ and $30$, respectively. The ideals are
given by
$$\left\{\aligned
  \mathfrak{a}_1 &=(\Phi_4, \Phi_8),\\
  \mathfrak{a}_2 &=(\Phi_4, \Phi_8, \Phi_{10}, \Phi_{4, 0}, \Phi_{1, 2}),\\
  \mathfrak{a}_3 &=(\Phi_4, \Phi_8, \Phi_{10}, \Phi_{3, 0}, \Phi_{0, 2},
                    \Phi_{14}, \Phi_{4, 0}, \Phi_{1, 2}, \Phi_{1, 3},
                    \Phi_{4, 1}),
\endaligned\right.\eqno{(1.16)}$$
and the parameter space $\{ (\lambda, \mu, \gamma) \} \cong \mathbb{C}^4$.
They form a variation of the $E_8$-singularity structure over the modular
curve $X$.}

\textbf{Theorem 1.2.} (Main Theorem 2) (A variation of the $E_8$-singularity
structure over the modular curve $X$: geometric version) {\it There is a
morphism of schemes
$$f: C_Y/G \rightarrow \text{Spec} \left(\mathbb{C}[\Phi_{12}, \Phi_{20},
  \Phi_{30}]/(\Phi_{20}^3-\Phi_{30}^2-1728 \Phi_{12}^5)\right) \eqno{(1.17)}$$
over the modular curve $X$. In particular, there are infinitely many such
triples $(\Phi_{12}, \Phi_{20}, \Phi_{30})$ $=(\Phi_{12}^{\lambda},
\Phi_{20}^{\mu}, \Phi_{30}^{\gamma})$ whose parameter space $\{ (\lambda,
\mu, \gamma) \}$ $\cong \mathbb{C}^4$. They form a variation of the
$E_8$-singularity structure over the modular curve $X$.}

  Theorem 1.1 and Theorem 1.2 show that there exist infinitely many kinds of
distinct constructions of the $E_8$-singularity: one and only one is given by
the Kleinian singularity $\mathbb{C}^2/\text{SL}(2, 5)$ (see \cite{K}), i.e.,
the icosahedral singularity, the other infinitely many kinds of constructions
are given from the quotient $C_Y/\text{SL}(2, 13)$ over the modular curve $X$.
Hence, the equation of $E_8$-singularity possesses infinitely many kinds of
distinct modular parametrizations.

  In his talk at ICM 1970 \cite{Br2}, Brieskorn showed how to construct
the singularity of type $ADE$ directly from the simple complex Lie group
of the same type. Namely, assume that $G$ is of type $ADE$, Brieskorn
proved a conjecture made by Grothendieck that the intersection of a
transversal slice to the sub-regular unipotent orbit with the unipotent
variety has a simple surface singularity of the same type as $G$. A
fuller treatment was given by Slodowy (see \cite{Br2} and \cite{Sl1}). A
clarification of the occurrence of the polyhedral groups in Brieskorn's
construction (see \cite{GrP} and \cite{GrP2}),  and thus a direct
relationship between the simple Lie groups and the finite subgroups of
$\text{SL}(2, \mathbb{C})$, was achieved by Kronheimer (see \cite{Kr1} and
\cite{Kr2}) using differential geometric methods. His construction starts
directly from the finite subgroups of $\text{SL}(2, \mathbb{C})$ and uses
hyper-K\"{a}hler quotient constructions. Kronheimer also gave an algebraic
approach using McKay correspondence. However, Brieskorn had still written
at the end of \cite{Br2}: ``Thus we see that there is a relation between
exotic spheres, the icosahedron and $E_8$. But I still do not understand
why the regular polyhedra come in.'' (see also \cite{Gr}, \cite{GrP},
\cite{GrP2} and \cite{Br3}). On the other hand, Arnol'd pointed out that
the theory of singularities is even linked (in a quite mysterious way) to
the classification of regular polyhedra in three-dimensional Euclidean
space (see \cite{Ar3}, p. 43). In his survey article on Platonic solids,
Kleinian singularities and Lie groups \cite{Sl2}, Slodowy found that the
objects of these different classifications are related to each other by
mathematical constructions. However, up to now, these constructions do
not explain why the different classifications should be related at all.

  As a consequence, Theorem 1.1 and Theorem 1.2 show that the
$E_8$-singularity is not necessarily the Kleinian icosahedral singularity.
That is, the icosahedron does not necessarily appear in the triple (exotic
spheres, icosahedron, $E_8$) of Brieskorn \cite{Br2}. The group
$\text{SL}(2, 13)$ can take its place and there are infinitely many kinds
of the other triples (exotic spheres, $\text{SL}(2, 13)$, $E_8$). The link
of these infinitely many kinds of distinct constructions of the $E_8$-singularity:
$\mathbb{C}^2/\text{SL}(2, 5)$ and a variation of the $E_8$-singularity
structure
$$C_Y/G \rightarrow \text{Spec}\left(\mathbb{C}[\Phi_{12}, \Phi_{20},
  \Phi_{30}]/(\Phi_{20}^3-\Phi_{30}^2-1728 \Phi_{12}^5)\right)$$
over the modular curve $X$ gives the same Poincar\'{e} homology $3$-sphere,
whose higher dimensional lifting:
$$z_1^5+z_2^3+z_3^2+z_4^2+z_5^2=0, \quad \sum_{i=1}^{5} z_i \overline{z_i}=1,
  \quad z_i \in \mathbb{C} \quad (1 \leq i \leq 5) \eqno{(1.18)}$$
gives the Milnor's standard generator of $\Theta_7$ (which is the version of
differential topology). Hence, this gives a negative answer to Arnol'd and
Brieskorn's questions about the mysterious relation between the icosahedron
and $E_8$, and the relation between Platonic solids, Kleinian singularities
and Lie groups appearing in Slodowy's survey \cite{Sl2} can be replaced by
the relation between $\text{SL}(2, 13)$, a variation of the $E_8$-singularity
structure
$$C_Y/\text{SL}(2, 13) \rightarrow \text{Spec} \left(\mathbb{C}[\Phi_{12},
 \Phi_{20}, \Phi_{30}]/(\Phi_{20}^3-\Phi_{30}^2-1728 \Phi_{12}^5) \right)$$
over the modular curve $X$ and $E_8$.

  Moreover, Theorem 1.1 and Theorem 1.2 can be extended to the following
two kinds of singularities:
$$\left\{\aligned
  &Q_{18}: & x^8+y^3+xz^2=0,\\
  &E_{20}: & x^{11}+y^3+z^2=0,
\endaligned\right.\eqno{(1.19)}$$
where $Q_{18}$ and $E_{20}$ are two bimodal singularities in the
pyramids of $14$ exceptional singularities (see \cite{Ar2}, p.255).
Theorem 1.1, Theorem 1.2, as well as Theorem 1.3 and Theorem 1.4
show that three different algebraic surfaces, the equations of $E_8$,
$Q_{18}$ and $E_{20}$-singularities can be realized from the same
quotients $C_Y/\text{SL}(2, 13)$ over the modular curve $X$ and have
the same modular parametrizations.

\textbf{Theorem 1.3.} (Variations of $Q_{18}$ and $E_{20}$-singularity
structures over the modular curve $X$: algebraic version) {\it The equations
of $Q_{18}$ and $E_{20}$-singularities
$$\Phi_{32}^3-\Phi_{12} \Phi_{42}^2-1728 \Phi_{12}^8=0, \quad
  \Phi_{44}^3-\Phi_{12}^4 \Phi_{42}^2-1728 \Phi_{12}^{11}=0$$
possess an infinitely many kinds of distinct modular parametrizations
$($with the cardinality of the continuum in ZFC set theory$)$
$$(\Phi_{12}, \Phi_{32}, \Phi_{42}, \Phi_{44})=(\Phi_{12}^{\lambda},
   \Phi_{32}^{\mu}, \Phi_{42}^{\gamma}, \Phi_{44})\eqno{(1.20)}$$
over the modular curve $X$ as follows$:$
$$\left\{\aligned
  \Phi_{12}^{\lambda} &=\lambda \Phi_{3, 0}+(1-\lambda) \Phi_{0, 2}
              \quad \text{mod $\mathfrak{a}_1$},\\
  \Phi_{32}^{\mu} &=\mu_1 \Phi_{8, 0}+\mu_2 \Phi_{5, 2}+(1-\mu_1-\mu_2)
              \Phi_{2, 4} \quad \text{mod $\mathfrak{a}_3$},\\
  \Phi_{42}^{\gamma} &=\gamma_1 \Phi_{0, 7}+\gamma_2 \Phi_{3, 5}+\gamma_3
              \Phi_{6, 3}+(1-\gamma_1-\gamma_2-\gamma_3) \Phi_{9, 1}
              \quad \text{mod $\mathfrak{a}_4$},\\
  \Phi_{44} &=\Phi_{11, 0} \quad \text{mod $\mathfrak{a}_4$},
\endaligned\right.\eqno{(1.21)}$$
where $\Phi_{12}$, $\Phi_{32}$, $\Phi_{42}$ and $\Phi_{44}$ are invariant
homogeneous polynomials of degree $12$, $32$, $42$ and $44$, respectively.
The ideals are given by
$$\left\{\aligned
  \mathfrak{a}_1 &=(\Phi_4, \Phi_8),\\
  \mathfrak{a}_3 &=(\Phi_4, \Phi_8, \Phi_{10}, \Phi_{3, 0}, \Phi_{0, 2},
  \Phi_{14}, \Phi_{4, 0}, \Phi_{1, 2}, \Phi_{1, 3}, \Phi_{4, 1}),\\
  \mathfrak{a}_4 &=(\Phi_4, \Phi_8, \Phi_{10}, \Phi_{3, 0}, \Phi_{0, 2},
  \Phi_{14}, \Phi_{4, 0}, \Phi_{1, 2}, \Phi_{1, 3}, \Phi_{4, 1}, \Phi_{1, 5},
  \Phi_{4, 3}, \Phi_{7, 1}),
\endaligned\right.\eqno{(1.22)}$$
and the parameter space $\{ (\lambda, \mu, \gamma) \} \cong \mathbb{C}^6$.
They form variations of $Q_{18}$ and $E_{20}$-singularity structures
over the modular curve $X$.}

\textbf{Theorem 1.4.} (Variations of $Q_{18}$ and $E_{20}$-singularity
structures over the modular curve $X$: geometric version) {\it There are two
morphisms from the cone $C_Y$ over $Y$ to the $Q_{18}$ and $E_{20}$-singularities:
$$f_1: C_Y/G \rightarrow \text{Spec} \left(\mathbb{C}[\Phi_{12}, \Phi_{32},
  \Phi_{42}]/(\Phi_{32}^3-\Phi_{12} \Phi_{42}^2-1728 \Phi_{12}^8)\right)
  \eqno{(1.23)}$$
and
$$f_2: C_Y/G \rightarrow \text{Spec} \left(\mathbb{C}[\Phi_{12}, \Phi_{42},
  \Phi_{44}]/(\Phi_{44}^3-\Phi_{12}^4 \Phi_{42}^2-1728 \Phi_{12}^{11})\right)
  \eqno{(1.24)}$$
over the modular curve $X$. In particular, there are infinitely many such
triples $(\Phi_{12}, \Phi_{32}, \Phi_{42})$ $=$ $(\Phi_{12}^{\lambda}, \Phi_{32}^{\mu},
\Phi_{42}^{\gamma})$ whose parameter space $\{ (\lambda, \mu, \gamma) \}$
$\cong \mathbb{C}^6$. They form variations of $Q_{18}$ and $E_{20}$-singularity
structures over the modular curve $X$.}

  In fact, Klein had noticed the similarity between the relation of
the equation $x^5+y^3+z^2=0$ to the icosahedral group $\text{PSL}(2, 5)$
and the relation of the equation $x^7+y^3+z^2=0$  to the group
$\text{PSL}(2, 7)$ (see \cite{K1}, \cite{K2}, \cite{KF1} and \cite{KF2}).
This is the starting point of the work of Dolgachev (see \cite{Do}) to
which Arnol'd was referring when he spoke about the wonderful coincidences
with Lobatchevsky triangles and automorphic functions (see \cite{Ar1}).
The normal form of Arnol'd for the quasi-homogeneous singularity $E_{12}$
in three variables is $x^7+y^3+z^2$, which can be realized as the quotient
conical singularity as follows (see \cite{BrPR} and \cite{Do}): The
canonical model $Y$ of the modular curve $X(7)$ in $\mathbb{CP}^2$ is
the Klein quartic given by the homogeneous equation
$z_1^3 z_2+z_2^3 z_3+z_3^3 z_1=0$. The finite group $\text{PSL}(2, 7)$
acts linearly on $\mathbb{C}^3$ and on $\mathbb{CP}^2$ leaving invariant
$Y \subset \mathbb{CP}^2$ and the cone $C_Y \subset \mathbb{C}^3$.
Calculations of invariants by Klein and Gordan imply:
$$\left[\mathbb{C}[z_1, z_2, z_3]/(z_1^3 z_2+z_2^3 z_3+z_3^3 z_1)
  \right]^{\text{PSL}(2, 7)} \cong \mathbb{C}[x, y, z]/(x^7+y^3+z^2).
  \eqno{(1.25)}$$
This algebraic result can be interpreted geometrically as follows:
The affine algebraic surface defined by the equation $x^7+y^3+z^2=0$
is the quotient of the cone $C_Y$ by the group $\text{PSL}(2, 7)$
over the modular curve $X(7)$, where $C_Y$ is the cone over $Y$.
Similarly, Klein also obtained the structure of the $\mathbb{C}$-algebra
of $\mathbb{C}[z_1, z_2]^{\text{SL}(2, 5)}$ of $\text{SL}(2, 5)$-invariant
polynomials on $\mathbb{C}^2$:
$$\mathbb{C}[z_1, z_2]^{\text{SL}(2, 5)} \cong \mathbb{C}[x, y, z]/
  (x^5+y^3+z^2).\eqno{(1.26)}$$
This algebraic result can also be interpreted geometrically as follows:
The affine algebraic surface defined by the equation $x^5+y^3+z^2=0$ is
the quotient of the cone $C_Y$ by the group $\text{SL}(2, 5)$ over the
modular curve $X(5)$, where $Y=\mathbb{CP}^1$ is the canonical model
of the modular curve $X(5)$ and $C_Y$ is a cone over $Y$. Therefore,
(1.14), (1.15), (1.17), (1.20), (1.21), (1.23), (1.24), (1.25) and
(1.26) give a complete and unified description for the relation between
the $G$-invariant homogeneous polynomials and the associated singularities
corresponding to the genus zero modular curves $X_0(N)$, where $N=5$, $7$,
$13$ and $G=\text{SL}(2, 5)$, $\text{PSL}(2, 7)$ and $\text{SL}(2, 13)$,
respectively:
$$\begin{matrix}
  &C_Y/\text{SL}(2, 5) \quad & C_Y/\text{PSL}(2, 7) \quad & C_Y/\text{SL}(2, 13)\\
  &\downarrow          \quad &\downarrow            \quad & \downarrow\\
  &X(5)                \quad &X(7)                  \quad & X(13)\\
  &E_8\text{-singularity} \quad &E_{12}\text{-singularity} \quad
  &\text{$E_8$, $Q_{18}$ and $E_{20}$-singularities}
\end{matrix}\eqno{(1.27)}$$
Here, $Y=\mathbb{CP}^1$, Klein quartic curve and our curve Y given by
(1.12), respectively.

  Finally, recall that there is a decomposition formula of the
elliptic modular function $j$ in terms of the icosahedral invariants
$f$, $H$ and $T$ of degrees $12$, $20$ and $30$ over the modular
curve $X(5)$ (see section two, in particular (2.8) for the details):
$$\aligned
 &j(z): j(z)-1728: 1\\
=&H(x_1(z), x_2(z))^3: -T(x_1(z), x_2(z))^2: f(x_1(z), x_2(z))^5,
\endaligned\eqno{(1.28)}$$
which was discovered by Klein (see \cite{K}, \cite{KF1} and \cite{KF2})
and later by Ramanujan (see \cite{Du}). In contrast with (1.28), we
have the following $j$-function version:

\textbf{Theorem 1.5.} (Main Theorem 3) (A variation of the structure of
decomposition formulas of the elliptic modular functions $j$ over the
modular curve $X$) {\it There are infinitely many kinds of distinct
decomposition formulas of the elliptic modular function $j$ in terms
of the invariants $\Phi_{12}$, $\Phi_{20}$ and $\Phi_{30}$ over the
modular curve $X$:
$$j(z): j(z)-1728: 1=\Phi_{20}^3: \Phi_{30}^2: \Phi_{12}^5,
  \eqno{(1.29)}$$
where $(\Phi_{12}, \Phi_{20}, \Phi_{30})=(\Phi_{12}^{\lambda},
\Phi_{20}^{\mu}, \Phi_{30}^{\gamma})$ are given by (1.15). They form a
variation of the structure of decomposition formulas of the elliptic
modular functions $j$ over the modular curve $X$.}

  In fact, these infinitely many kinds of distinct decompositions (1.29)
and (1.28) have the same form, i.e., the degrees of the invariant polynomials
are $12$, $20$ and $30$, respectively. However, they have the different
geometric interpretation: one and only one is over the modular curve $X(5)$,
the other infinitely many kinds of decompositions (which form a variation of
the structure of decomposition formulas) are over the modular curve $X(13)$.
They also have the different algebraic interpretation: one and only one is
invariant under the group $\text{SL}(2, 5)$, the other infinitely many kinds of
decompositions (which form a variation of the structure of decomposition formulas)
are invariant under the group $\text{SL}(2, 13)$.

  This paper consists of five sections. In section two, we revisit
the standard construction of the $E_8$-singularity as the well-known
Kleinian icosahedral singularity. In section three, we study the
invariant theory and modular forms for $\text{SL}(2, 13)$. In
particular, we construct a system of invariants $\Phi_{m, n}$
for $\text{SL}(2, 13)$. These invariants are modular forms over the
modular curve $X$. In section four, we find three invariant homogeneous
polynomials $\Phi_{12}$, $\Phi_{20}$ and $\Phi_{30}$ among those modular
forms. They satisfy the equation of $E_8$-singularity. Thus we obtain the
ring homomorphism from
$\mathbb{C}[\Phi_{12}, \Phi_{20}, \Phi_{30}]/(\Phi_{20}^3-\Phi_{30}^2-1728
\Phi_{12}^5)$ to the ring of invariants
$\left[\mathbb{C}[z_1, z_2, z_3, z_4, z_5, z_6]/I\right]^{\text{SL}(2, 13)}$,
where $I$ is an ideal generated by a system of algebraic relations satisfied
by the variables $z_1$, $\ldots$, $z_6$ and $Y$ is the algebraic curve
corresponding to the ideal $I$. In particular, there are infinitely many
such triples $(\Phi_{12}^{\lambda}, \Phi_{20}^{\mu}, \Phi_{30}^{\gamma})$,
which form a variation of the $E_8$-singularity structure over the modular
curve $X$. We give its algebraic version, geometric version, $j$-function
version and the version of differential topology: Poincar\'{e} homology
$3$-sphere as well as its higher dimensional lifting, i.e., Milnor's exotic
$7$-sphere. Therefore, we give a different construction of the $E_8$-singularity
coming from a quotient $C_Y/\text{SL}(2, 13)$ over the modular curve $X$, where
$C_Y$ is the cone over the algebraic curve $Y$. In section five, we extend our
work to the cases of $Q_{18}$ and $E_{20}$-singularities and obtain variations
of $Q_{18}$ and $E_{20}$-singularity structures over the modular curve $X$.

\textbf{Acknowledgements}. The author would like to thank Pierre Deligne
for his very detailed and helpful comments as well as his patience.

\begin{center}
{\large\bf 2. Standard construction: $E_8$-singularity as the
              icosahedral singularity}
\end{center}

  Let us recall some classical result on the relation between the
icosahedron and the $E_8$-singularity (see \cite{Mc}). Starting
with the polynomial invariants of the finite subgroup of
$\text{SL}(2, \mathbb{C})$, a surface is defined from the single
syzygy which relates the three polynomials in two variables. This
surface has a singularity at the origin; the singularity can be
resolved by constructing a smooth surface which is isomorphic to
the original one except for a set of component curves which form
the pre-image of the origin. The components form a Dynkin curve
and the matrix of their intersections is the negative of the Cartan
matrix for the appropriate Lie algebra. The Dynkin curve is the
dual of the Dynkin graph. For example, if $\Gamma$ is the binary
icosahedral group, the corresponding Dynkin curve is that of $E_8$,
and $\mathbb{C}^2/\Gamma \subset \mathbb{C}^3$ is the set of zeros
of the equation
$$x^5+y^3+z^2=0.\eqno{(2.1)}$$
The link of this $E_8$-singularity, the Poincar\'{e} homology
$3$-sphere (see \cite{KS}), has a higher dimensional lifting:
$$z_1^5+z_2^3+z_3^2+z_4^2+z_5^2=0, \quad \sum_{i=1}^{5} z_i \overline{z_i}=1,
  \quad z_i \in \mathbb{C} \quad (1 \leq i \leq 5),\eqno{(2.2)}$$
which is the Brieskorn description of one of Milnor's exotic
$7$-dimensional spheres. In fact, it is an exotic $7$-sphere
representing Milnor's standard generator of $\Theta_7$ (see
\cite{Br1}, \cite{Br2} and \cite{Hi}).

  In his celebrated book \cite{K}, Klein gave a parametric solution
of the above singularity (2.1) by homogeneous polynomials $T$, $H$,
$f$ in two variables of degrees $30$, $20$, $12$ with integral
coefficients, where
$$f=z_1 z_2 (z_1^{10}+11 z_1^5 z_2^5-z_2^{10}),$$
$$H=\frac{1}{121} \begin{vmatrix}
    \frac{\partial^2 f}{\partial z_1^2} &
    \frac{\partial^2 f}{\partial z_1 \partial z_2}\\
    \frac{\partial^2 f}{\partial z_2 \partial z_1} &
    \frac{\partial^2 f}{\partial z_2^2}
    \end{vmatrix}
  =-(z_1^{20}+z_2^{20})+228 (z_1^{15} z_2^5-z_1^5 z_2^{15})
   -494 z_1^{10} z_2^{10},$$
$$T=-\frac{1}{20} \begin{vmatrix}
    \frac{\partial f}{\partial z_1} &
    \frac{\partial f}{\partial z_2}\\
    \frac{\partial H}{\partial z_1} &
    \frac{\partial H}{\partial z_2}
    \end{vmatrix}
  =(z_1^{30}+z_2^{30})+522 (z_1^{25} z_2^5-z_1^5 z_2^{25})
   -10005 (z_1^{20} z_2^{10}+z_1^{10} z_2^{20}).$$
They satisfy the famous (binary) icosahedral equation
$$T^2+H^3=1728 f^5.\eqno{(2.3)}$$
In fact, $f$, $H$ and $T$ are invariant polynomials under the
action of the binary icosahedral group. The above equation (2.3)
is closely related to Hermite's celebrated work (see \cite{He})
on the resolution of the quintic equations by elliptic modular
functions of order five. Essentially the same relation had been
found a few years earlier by Schwarz (see \cite{Sch}), who
considered three polynomials $\varphi_{12}$, $\varphi_{20}$ and
$\varphi_{30}$ whose roots correspond to the vertices, the
midpoints of the faces and the midpoints of the edges of an
icosahedron inscribed in the Riemann sphere. He obtained the
identity $\varphi_{20}^3-1728 \varphi_{12}^5=\varphi_{30}^2$.
We see this identity as well as (2.3) as the defining relation
between three generators $f$, $H$ and $T$ of the ring of invariants
$\mathbb{C}[z_1, z_2]^{\Gamma}$ of the binary icosahedral group
$\Gamma$ acting on $\mathbb{C}^2$, and we identify this ring with
the ring of functions on the affine variety $\mathbb{C}^2/\Gamma$
embedded in $\mathbb{C}^3$ and given by such an equation (see
\cite{Br3}). Namely,
$$\mathbb{C}[z_1, z_2]^{\Gamma} \cong \mathbb{C}[f, H, T]/
  (T^2+H^3-1728 f^5).\eqno{(2.4)}$$
Thus we see that from the very beginning there was a close
relation between the $E_8$-singularity and the icosahedron.
Moreover, the icosahedral equation (2.3) can be interpreted
in terms of modular forms which was also known by Klein (see
\cite{KF1}, p. 631). Let $x_1(z)=\eta(z) a(z)$ and
$x_2(z)=\eta(z) b(z)$, where
$$a(z)=e^{-\frac{3 \pi i}{10}} \theta \begin{bmatrix}
       \frac{3}{5}\\ 1 \end{bmatrix}(0, 5z), \quad
  b(z)=e^{-\frac{\pi i}{10}} \theta \begin{bmatrix}
       \frac{1}{5}\\ 1 \end{bmatrix}(0, 5z)$$
are theta constants of order five and
$\eta(z):=q^{\frac{1}{24}} \prod_{n=1}^{\infty} (1-q^n)$ with
$q=e^{2 \pi i z}$ is the Dedekind eta function which are all
defined in the upper-half plane
$\mathbb{H}=\{ z \in \mathbb{C}: \text{Im}(z)>0 \}$. Then
$$\left\{\aligned
  f(x_1(z), x_2(z)) &=-\Delta(z),\\
  H(x_1(z), x_2(z)) &=-\eta(z)^8 \Delta(z)E_4(z),\\
  T(x_1(z), x_2(z)) &=\Delta(z)^2 E_6(z),
\endaligned\right.\eqno{(2.5)}$$
where
$$E_4(z):=\frac{1}{2} \sum_{m, n \in \mathbb{Z}, (m, n)=1}
          \frac{1}{(mz+n)^4}, \quad
  E_6(z):=\frac{1}{2} \sum_{m, n \in \mathbb{Z}, (m, n)=1}
          \frac{1}{(mz+n)^6}$$
are Eisenstein series of weight $4$ and $6$, and
$\Delta(z)=\eta(z)^{24}$ is the discriminant. The relations
$$j(z)=\frac{E_4(z)^3}{\Delta(z)}=\frac{H(x_1(z), x_2(z))^3}
        {f(x_1(z), x_2(z))^5},\eqno{(2.6)}$$
$$j(z)-1728=\frac{E_6(z)^2}{\Delta(z)}=-\frac{T(x_1(z), x_2(z))^2}
            {f(x_1(z), x_2(z))^5}\eqno{(2.7)}$$
give the icosahedral equation (2.3) in terms of theta constants
of order five. Hence, we have the following decomposition formula
of the elliptic modular function $j$ in terms of the icosahedral
invariants $f$, $H$ and $T$ over the modular curve $X(5)$:
$$\aligned
 &j(z): j(z)-1728: 1\\
=&H(x_1(z), x_2(z))^3: -T(x_1(z), x_2(z))^2: f(x_1(z), x_2(z))^5.
\endaligned\eqno{(2.8)}$$

\begin{center}
{\large\bf 3. Invariant theory and modular forms for $\text{SL}(2, 13)$}
\end{center}

  The representation of $\text{SL}(2, 13)$ which we will consider is
the unique six-dimensional irreducible complex representation for which
the eigenvalues of $\left(\begin{matrix} 1 & 1\\ 0 & 1 \end{matrix}\right)$
are the $\exp (\underline{a} . 2 \pi i/13)$ for $\underline{a}$ a
non-square mod $13$. We will give an explicit realization of this
representation. This explicit realization will play a major role for
giving a complete system of invariants associated to $\text{SL}(2, 13)$.
At first, we will study the six-dimensional representation of
the finite group $\text{SL}(2, 13)$ of order $2184$, which acts
on the five-dimensional projective space
$\mathbb{P}^5=\{ (z_1, z_2, z_3, z_4, z_5, z_6): z_i \in \mathbb{C}
 \quad (i=1, 2, 3, 4, 5, 6) \}$. This representation is defined
over the cyclotomic field $\mathbb{Q}(e^{\frac{2 \pi i}{13}})$.
Put
$$S=-\frac{1}{\sqrt{13}} \begin{pmatrix}
  \zeta^{12}-\zeta & \zeta^{10}-\zeta^3 & \zeta^4-\zeta^9 &
  \zeta^5-\zeta^8 & \zeta^2-\zeta^{11} & \zeta^6-\zeta^7\\
  \zeta^{10}-\zeta^3 & \zeta^4-\zeta^9 & \zeta^{12}-\zeta &
  \zeta^2-\zeta^{11} & \zeta^6-\zeta^7 & \zeta^5-\zeta^8\\
  \zeta^4-\zeta^9 & \zeta^{12}-\zeta & \zeta^{10}-\zeta^3 &
  \zeta^6-\zeta^7 & \zeta^5-\zeta^8 & \zeta^2-\zeta^{11}\\
  \zeta^5-\zeta^8 & \zeta^2-\zeta^{11} & \zeta^6-\zeta^7 &
  \zeta-\zeta^{12} & \zeta^3-\zeta^{10} & \zeta^9-\zeta^4\\
  \zeta^2-\zeta^{11} & \zeta^6-\zeta^7 & \zeta^5-\zeta^8 &
  \zeta^3-\zeta^{10} & \zeta^9-\zeta^4 & \zeta-\zeta^{12}\\
  \zeta^6-\zeta^7 & \zeta^5-\zeta^8 & \zeta^2-\zeta^{11} &
  \zeta^9-\zeta^4 & \zeta-\zeta^{12} & \zeta^3-\zeta^{10}
\end{pmatrix}\eqno{(3.1)}$$
and
$$T=\text{diag}(\zeta^7, \zeta^{11}, \zeta^8, \zeta^6,
    \zeta^2, \zeta^5),\eqno{(3.2)}$$
where $\zeta=\exp(2 \pi i/13)$. We have
$$S^2=-I, \quad T^{13}=(ST)^3=I.\eqno{(3.3)}$$
In \cite{Y1}, we put $P=S T^{-1} S$ and $Q=S T^3$. Then
$(Q^3 P^4)^3=-I$ (see \cite{Y1}, the proof of Theorem 3.1).
Let $G=\langle S, T \rangle$, then $G \cong \text{SL}(2, 13)$.

  Put $\theta_1=\zeta+\zeta^3+\zeta^9$,
$\theta_2=\zeta^2+\zeta^6+\zeta^5$,
$\theta_3=\zeta^4+\zeta^{12}+\zeta^{10}$,
and $\theta_4=\zeta^8+\zeta^{11}+\zeta^7$. We find that
$$\left\{\aligned
  &\theta_1+\theta_2+\theta_3+\theta_4=-1,\\
  &\theta_1 \theta_2+\theta_1 \theta_3+\theta_1 \theta_4+
   \theta_2 \theta_3+\theta_2 \theta_4+\theta_3 \theta_4=2,\\
  &\theta_1 \theta_2 \theta_3+\theta_1 \theta_2 \theta_4+
   \theta_1 \theta_3 \theta_4+\theta_2 \theta_3 \theta_4=4,\\
  &\theta_1 \theta_2 \theta_3 \theta_4=3.
\endaligned\right.$$
Hence, $\theta_1$, $\theta_2$, $\theta_3$ and $\theta_4$ satisfy
the quartic equation $z^4+z^3+2 z^2-4z+3=0$,
which can be decomposed as two quadratic equations
$$\left(z^2+\frac{1+\sqrt{13}}{2} z+\frac{5+\sqrt{13}}{2}\right)
  \left(z^2+\frac{1-\sqrt{13}}{2} z+\frac{5-\sqrt{13}}{2}\right)=0$$
over the real quadratic field $\mathbb{Q}(\sqrt{13})$. Therefore, the
four roots are given as follows:
$$\left\{\aligned
  \theta_1=\frac{1}{4} \left(-1+\sqrt{13}+\sqrt{-26+6 \sqrt{13}}\right),\\
  \theta_2=\frac{1}{4} \left(-1-\sqrt{13}+\sqrt{-26-6 \sqrt{13}}\right),\\
  \theta_3=\frac{1}{4} \left(-1+\sqrt{13}-\sqrt{-26+6 \sqrt{13}}\right),\\
  \theta_4=\frac{1}{4} \left(-1-\sqrt{13}-\sqrt{-26-6 \sqrt{13}}\right).
\endaligned\right.$$
Moreover, we find that
$$\left\{\aligned
  \theta_1+\theta_3+\theta_2+\theta_4 &=-1,\\
  \theta_1+\theta_3-\theta_2-\theta_4 &=\sqrt{13},\\
  \theta_1-\theta_3-\theta_2+\theta_4 &=-\sqrt{-13+2 \sqrt{13}},\\
  \theta_1-\theta_3+\theta_2-\theta_4 &=\sqrt{-13-2 \sqrt{13}}.
\endaligned\right.$$

  Let us study the action of $S T^{\nu}$ on $\mathbb{P}^5$, where
$\nu=0, 1, \ldots, 12$. Put
$$\alpha=\zeta+\zeta^{12}-\zeta^5-\zeta^8, \quad
   \beta=\zeta^3+\zeta^{10}-\zeta^2-\zeta^{11}, \quad
   \gamma=\zeta^9+\zeta^4-\zeta^6-\zeta^7.$$
We find that
$$\aligned
  &13 ST^{\nu}(z_1) \cdot ST^{\nu}(z_4)\\
=&\beta z_1 z_4+\gamma z_2 z_5+\alpha z_3 z_6+\\
 &+\gamma \zeta^{\nu} z_1^2+\alpha \zeta^{9 \nu} z_2^2+\beta \zeta^{3 \nu} z_3^2
  -\gamma \zeta^{12 \nu} z_4^2-\alpha \zeta^{4 \nu} z_5^2-\beta \zeta^{10 \nu}
  z_6^2+\\
 &+(\alpha-\beta) \zeta^{5 \nu} z_1 z_2+(\beta-\gamma) \zeta^{6 \nu} z_2 z_3
  +(\gamma-\alpha) \zeta^{2 \nu} z_1 z_3+\\
 &+(\beta-\alpha) \zeta^{8 \nu} z_4 z_5+(\gamma-\beta) \zeta^{7 \nu} z_5 z_6
  +(\alpha-\gamma) \zeta^{11 \nu} z_4 z_6+\\
 &-(\alpha+\beta) \zeta^{\nu} z_3 z_4-(\beta+\gamma) \zeta^{9 \nu} z_1 z_5
  -(\gamma+\alpha) \zeta^{3 \nu} z_2 z_6+\\
 &-(\alpha+\beta) \zeta^{12 \nu} z_1 z_6-(\beta+\gamma) \zeta^{4 \nu} z_2 z_4
  -(\gamma+\alpha) \zeta^{10 \nu} z_3 z_5.
\endaligned$$
$$\aligned
  &13 ST^{\nu}(z_2) \cdot ST^{\nu}(z_5)\\
=&\gamma z_1 z_4+\alpha z_2 z_5+\beta z_3 z_6+\\
 &+\alpha \zeta^{\nu} z_1^2+\beta \zeta^{9 \nu} z_2^2+\gamma \zeta^{3 \nu} z_3^2
  -\alpha \zeta^{12 \nu} z_4^2-\beta \zeta^{4 \nu} z_5^2-\gamma \zeta^{10 \nu}
  z_6^2+\\
 &+(\beta-\gamma) \zeta^{5 \nu} z_1 z_2+(\gamma-\alpha) \zeta^{6 \nu} z_2 z_3
  +(\alpha-\beta) \zeta^{2 \nu} z_1 z_3+\\
 &+(\gamma-\beta) \zeta^{8 \nu} z_4 z_5+(\alpha-\gamma) \zeta^{7 \nu} z_5 z_6
  +(\beta-\alpha) \zeta^{11 \nu} z_4 z_6+\\
 &-(\beta+\gamma) \zeta^{\nu} z_3 z_4-(\gamma+\alpha) \zeta^{9 \nu} z_1 z_5
  -(\alpha+\beta) \zeta^{3 \nu} z_2 z_6+\\
 &-(\beta+\gamma) \zeta^{12 \nu} z_1 z_6-(\gamma+\alpha) \zeta^{4 \nu} z_2 z_4
  -(\alpha+\beta) \zeta^{10 \nu} z_3 z_5.
\endaligned$$
$$\aligned
  &13 ST^{\nu}(z_3) \cdot ST^{\nu}(z_6)\\
=&\alpha z_1 z_4+\beta z_2 z_5+\gamma z_3 z_6+\\
 &+\beta \zeta^{\nu} z_1^2+\gamma \zeta^{9 \nu} z_2^2+\alpha \zeta^{3 \nu} z_3^2
  -\beta \zeta^{12 \nu} z_4^2-\gamma \zeta^{4 \nu} z_5^2-\alpha \zeta^{10 \nu}
  z_6^2+\\
 &+(\gamma-\alpha) \zeta^{5 \nu} z_1 z_2+(\alpha-\beta) \zeta^{6 \nu} z_2 z_3
  +(\beta-\gamma) \zeta^{2 \nu} z_1 z_3+\\
 &+(\alpha-\gamma) \zeta^{8 \nu} z_4 z_5+(\beta-\alpha) \zeta^{7 \nu} z_5 z_6
  +(\gamma-\beta) \zeta^{11 \nu} z_4 z_6+\\
 &-(\gamma+\alpha) \zeta^{\nu} z_3 z_4-(\alpha+\beta) \zeta^{9 \nu} z_1 z_5
  -(\beta+\gamma) \zeta^{3 \nu} z_2 z_6+\\
 &-(\gamma+\alpha) \zeta^{12 \nu} z_1 z_6-(\alpha+\beta) \zeta^{4 \nu} z_2 z_4
  -(\beta+\gamma) \zeta^{10 \nu} z_3 z_5.
\endaligned$$
Note that $\alpha+\beta+\gamma=\sqrt{13}$, we find that
$$\aligned
  &\sqrt{13} \left[ST^{\nu}(z_1) \cdot ST^{\nu}(z_4)+
   ST^{\nu}(z_2) \cdot ST^{\nu}(z_5)+ST^{\nu}(z_3) \cdot ST^{\nu}(z_6)\right]\\
 =&(z_1 z_4+z_2 z_5+z_3 z_6)+(\zeta^{\nu} z_1^2+\zeta^{9 \nu} z_2^2+\zeta^{3 \nu}
   z_3^2)-(\zeta^{12 \nu} z_4^2+\zeta^{4 \nu} z_5^2+\zeta^{10 \nu} z_6^2)+\\
  &-2(\zeta^{\nu} z_3 z_4+\zeta^{9 \nu} z_1 z_5+\zeta^{3 \nu} z_2 z_6)
   -2(\zeta^{12 \nu} z_1 z_6+\zeta^{4 \nu} z_2 z_4+\zeta^{10 \nu} z_3 z_5).
\endaligned$$
Let
$$\varphi_{\infty}(z_1, z_2, z_3, z_4, z_5, z_6)=\sqrt{13} (z_1 z_4+z_2 z_5+z_3 z_6)
  \eqno{(3.4)}$$
and
$$\varphi_{\nu}(z_1, z_2, z_3, z_4, z_5, z_6)=\varphi_{\infty}(ST^{\nu}(z_1, z_2,
                                              z_3, z_4, z_5, z_6))\eqno{(3.5)}$$
for $\nu=0, 1, \ldots, 12$. Then
$$\aligned
  \varphi_{\nu}
=&(z_1 z_4+z_2 z_5+z_3 z_6)+\zeta^{\nu} (z_1^2-2 z_3 z_4)+\zeta^{4 \nu}
  (-z_5^2-2 z_2 z_4)+\\
 &+\zeta^{9 \nu} (z_2^2-2 z_1 z_5)+\zeta^{3 \nu} (z_3^2-2 z_2 z_6)+
   \zeta^{12 \nu} (-z_4^2-2 z_1 z_6)+\\
 &+\zeta^{10 \nu} (-z_6^2-2 z_3 z_5).
\endaligned\eqno{(3.6)}$$
This leads us to define the following senary quadratic forms
(quadratic forms in six variables):
$$\left\{\aligned
  \mathbf{A}_0 &=z_1 z_4+z_2 z_5+z_3 z_6,\\
  \mathbf{A}_1 &=z_1^2-2 z_3 z_4,\\
  \mathbf{A}_2 &=-z_5^2-2 z_2 z_4,\\
  \mathbf{A}_3 &=z_2^2-2 z_1 z_5,\\
  \mathbf{A}_4 &=z_3^2-2 z_2 z_6,\\
  \mathbf{A}_5 &=-z_4^2-2 z_1 z_6,\\
  \mathbf{A}_6 &=-z_6^2-2 z_3 z_5.
\endaligned\right.\eqno{(3.7)}$$
Hence,
$$\sqrt{13} ST^{\nu}(\mathbf{A}_0)=\mathbf{A}_0+\zeta^{\nu} \mathbf{A}_1
  +\zeta^{4 \nu} \mathbf{A}_2+\zeta^{9 \nu} \mathbf{A}_3+\zeta^{3 \nu}
  \mathbf{A}_4+\zeta^{12 \nu} \mathbf{A}_5+\zeta^{10 \nu} \mathbf{A}_6.
  \eqno{(3.8)}$$
Let $H:=Q^5 P^2 \cdot P^2 Q^6 P^8 \cdot Q^5 P^2 \cdot P^3 Q$ where
$P=S T^{-1} S$ and $Q=S T^3$. Then (see \cite{Y2}, p.27)
$$H=\begin{pmatrix}
  0 &  0 &  0 & 0 & 0 & 1\\
  0 &  0 &  0 & 1 & 0 & 0\\
  0 &  0 &  0 & 0 & 1 & 0\\
  0 &  0 & -1 & 0 & 0 & 0\\
 -1 &  0 &  0 & 0 & 0 & 0\\
  0 & -1 &  0 & 0 & 0 & 0
\end{pmatrix}.\eqno{(3.9)}$$
Note that $H^6=1$ and $H^{-1} T H=-T^4$. Thus,
$\langle H, T \rangle \cong \mathbb{Z}_{13} \rtimes \mathbb{Z}_6$.
Hence, it is a maximal subgroup of order $78$ of $\text{PSL}(2, 13)$
with index $14$. We find that $\varphi_{\infty}^2$ is invariant under
the action of the maximal subgroup $\langle H, T \rangle$. Note that
$$\varphi_{\infty}=\sqrt{13} \mathbf{A}_0, \quad
  \varphi_{\nu}=\mathbf{A}_0+\zeta^{\nu} \mathbf{A}_1+
  \zeta^{4 \nu} \mathbf{A}_2+\zeta^{9 \nu} \mathbf{A}_3+
  \zeta^{3 \nu} \mathbf{A}_4+\zeta^{12 \nu} \mathbf{A}_5+
  \zeta^{10 \nu} \mathbf{A}_6$$
for $\nu=0, 1, \ldots, 12$. Let $w=\varphi^2$,
$w_{\infty}=\varphi_{\infty}^2$ and $w_{\nu}=\varphi_{\nu}^2$.
Then $w_{\infty}$, $w_{\nu}$ for $\nu=0, \ldots, 12$ form an
algebraic equation of degree fourteen, which is just the Jacobian
equation of degree fourteen (see \cite{K}, pp.161-162), whose
roots are these $w_{\nu}$ and $w_{\infty}$:
$$w^{14}+a_1 w^{13}+\cdots+a_{13} w+a_{14}=0.$$

  On the other hand, we have
$$\aligned
  &-13 \sqrt{13} ST^{\nu}(z_1) \cdot ST^{\nu}(z_2) \cdot ST^{\nu}(z_3)\\
 =&-r_4 (\zeta^{8 \nu} z_1^3+\zeta^{7 \nu} z_2^3+\zeta^{11 \nu} z_3^3)
   -r_2 (\zeta^{5 \nu} z_4^3+\zeta^{6 \nu} z_5^3+\zeta^{2 \nu} z_6^3)\\
  &-r_3 (\zeta^{12 \nu} z_1^2 z_2+\zeta^{4 \nu} z_2^2 z_3+\zeta^{10 \nu} z_3^2 z_1)
   -r_1 (\zeta^{\nu} z_4^2 z_5+\zeta^{9 \nu} z_5^2 z_6+\zeta^{3 \nu} z_6^2 z_4)\\
  &+2 r_1 (\zeta^{3 \nu} z_1 z_2^2+\zeta^{\nu} z_2 z_3^2+\zeta^{9 \nu} z_3 z_1^2)
   -2 r_3 (\zeta^{10 \nu} z_4 z_5^2+\zeta^{12 \nu} z_5 z_6^2+\zeta^{4 \nu} z_6 z_4^2)\\
  &+2 r_4 (\zeta^{7 \nu} z_1^2 z_4+\zeta^{11 \nu} z_2^2 z_5+\zeta^{8 \nu} z_3^2 z_6)
   -2 r_2 (\zeta^{6 \nu} z_1 z_4^2+\zeta^{2 \nu} z_2 z_5^2+\zeta^{5 \nu} z_3 z_6^2)+\\
  &+r_1 (\zeta^{3 \nu} z_1^2 z_5+\zeta^{\nu} z_2^2 z_6+\zeta^{9 \nu} z_3^2 z_4)
   +r_3 (\zeta^{10 \nu} z_2 z_4^2+\zeta^{12 \nu} z_3 z_5^2+\zeta^{4 \nu} z_1 z_6^2)+\\
  &+r_2 (\zeta^{6 \nu} z_1^2 z_6+\zeta^{2 \nu} z_2^2 z_4+\zeta^{5 \nu} z_3^2 z_5)
   +r_4 (\zeta^{7 \nu} z_3 z_4^2+\zeta^{11 \nu} z_1 z_5^2+\zeta^{8 \nu} z_2 z_6^2)+\\
  &+r_0 z_1 z_2 z_3+r_{\infty} z_4 z_5 z_6+\\
  &-r_4 (\zeta^{11 \nu} z_1 z_2 z_4+\zeta^{8 \nu} z_2 z_3 z_5+\zeta^{7 \nu} z_1 z_3 z_6)+\\
  &+r_2 (\zeta^{2 \nu} z_1 z_4 z_5+\zeta^{5 \nu} z_2 z_5 z_6+\zeta^{6 \nu} z_3 z_4 z_6)+\\
  &-3 r_4 (\zeta^{7 \nu} z_1 z_2 z_5+\zeta^{11 \nu} z_2 z_3 z_6+\zeta^{8 \nu} z_1 z_3 z_4)+\\
  &+3 r_2 (\zeta^{6 \nu} z_2 z_4 z_5+\zeta^{2 \nu} z_3 z_5 z_6+\zeta^{5 \nu} z_1 z_4 z_6)+\\
  &-r_3 (\zeta^{10 \nu} z_1 z_2 z_6+\zeta^{4 \nu} z_1 z_3 z_5+\zeta^{12 \nu} z_2 z_3 z_4)+\\
  &+r_1 (\zeta^{3 \nu} z_3 z_4 z_5+\zeta^{9 \nu} z_2 z_4 z_6+\zeta^{\nu} z_1 z_5 z_6),
\endaligned$$
where
$$r_0=2(\theta_1-\theta_3)-3(\theta_2-\theta_4), \quad
  r_{\infty}=2(\theta_4-\theta_2)-3(\theta_1-\theta_3),$$
$$r_1=\sqrt{-13-2 \sqrt{13}}, \quad r_2=\sqrt{\frac{-13+3 \sqrt{13}}{2}},$$
$$r_3=\sqrt{-13+2 \sqrt{13}}, \quad r_4=\sqrt{\frac{-13-3 \sqrt{13}}{2}}.$$
This leads us to define the following senary cubic forms (cubic forms in
six variables):
$$\left\{\aligned
  \mathbf{D}_0 &=z_1 z_2 z_3,\\
  \mathbf{D}_1 &=2 z_2 z_3^2+z_2^2 z_6-z_4^2 z_5+z_1 z_5 z_6,\\
  \mathbf{D}_2 &=-z_6^3+z_2^2 z_4-2 z_2 z_5^2+z_1 z_4 z_5+3 z_3 z_5 z_6,\\
  \mathbf{D}_3 &=2 z_1 z_2^2+z_1^2 z_5-z_4 z_6^2+z_3 z_4 z_5,\\
  \mathbf{D}_4 &=-z_2^2 z_3+z_1 z_6^2-2 z_4^2 z_6-z_1 z_3 z_5,\\
  \mathbf{D}_5 &=-z_4^3+z_3^2 z_5-2 z_3 z_6^2+z_2 z_5 z_6+3 z_1 z_4 z_6,\\
  \mathbf{D}_6 &=-z_5^3+z_1^2 z_6-2 z_1 z_4^2+z_3 z_4 z_6+3 z_2 z_4 z_5,\\
  \mathbf{D}_7 &=-z_2^3+z_3 z_4^2-z_1 z_3 z_6-3 z_1 z_2 z_5+2 z_1^2 z_4,\\
  \mathbf{D}_8 &=-z_1^3+z_2 z_6^2-z_2 z_3 z_5-3 z_1 z_3 z_4+2 z_3^2 z_6,\\
  \mathbf{D}_9 &=2 z_1^2 z_3+z_3^2 z_4-z_5^2 z_6+z_2 z_4 z_6,\\
  \mathbf{D}_{10} &=-z_1 z_3^2+z_2 z_4^2-2 z_4 z_5^2-z_1 z_2 z_6,\\
  \mathbf{D}_{11} &=-z_3^3+z_1 z_5^2-z_1 z_2 z_4-3 z_2 z_3 z_6+2 z_2^2 z_5,\\
  \mathbf{D}_{12} &=-z_1^2 z_2+z_3 z_5^2-2 z_5 z_6^2-z_2 z_3 z_4,\\
  \mathbf{D}_{\infty}&=z_4 z_5 z_6.
\endaligned\right.\eqno{(3.10)}$$
Then
$$\aligned
  &-13 \sqrt{13} ST^{\nu}(\mathbf{D}_0)\\
 =&r_0 \mathbf{D}_0+r_1 \zeta^{\nu} \mathbf{D}_1+
   r_2 \zeta^{2 \nu} \mathbf{D}_2+
   r_1 \zeta^{3 \nu} \mathbf{D}_3+r_3 \zeta^{4 \nu} \mathbf{D}_4+\\
  &+r_2 \zeta^{5 \nu} \mathbf{D}_5+r_2 \zeta^{6 \nu} \mathbf{D}_6+
   r_4 \zeta^{7 \nu} \mathbf{D}_7+r_4 \zeta^{8 \nu} \mathbf{D}_8+\\
  &+r_1 \zeta^{9 \nu} \mathbf{D}_9+r_3 \zeta^{10 \nu} \mathbf{D}_{10}
   +r_4 \zeta^{11 \nu} \mathbf{D}_{11}+r_3 \zeta^{12 \nu} \mathbf{D}_{12}
   +r_{\infty} \mathbf{D}_{\infty}.
\endaligned$$
$$\aligned
  &-13 \sqrt{13} ST^{\nu}(\mathbf{D}_{\infty})\\
 =&r_{\infty} \mathbf{D}_0-r_3 \zeta^{\nu} \mathbf{D}_1-
   r_4 \zeta^{2 \nu} \mathbf{D}_2-r_3 \zeta^{3 \nu} \mathbf{D}_3+
   r_1 \zeta^{4 \nu} \mathbf{D}_4+\\
  &-r_4 \zeta^{5 \nu} \mathbf{D}_5-r_4 \zeta^{6 \nu} \mathbf{D}_6+
   r_2 \zeta^{7 \nu} \mathbf{D}_7+r_2 \zeta^{8 \nu} \mathbf{D}_8+\\
  &-r_3 \zeta^{9 \nu} \mathbf{D}_9+r_1 \zeta^{10 \nu} \mathbf{D}_{10}+
   r_2 \zeta^{11 \nu} \mathbf{D}_{11}+r_1 \zeta^{12 \nu} \mathbf{D}_{12}-
   r_0 \mathbf{D}_{\infty}.
\endaligned$$

  Let
$$\delta_{\infty}(z_1, z_2, z_3, z_4, z_5, z_6)
 =13^2 (z_1^2 z_2^2 z_3^2+z_4^2 z_5^2 z_6^2)\eqno{(3.11)}$$
and
$$\delta_{\nu}(z_1, z_2, z_3, z_4, z_5, z_6)
 =\delta_{\infty}(ST^{\nu}(z_1, z_2, z_3, z_4, z_5, z_6))\eqno{(3.12)}$$
for $\nu=0, 1, \ldots, 12$. Then
$$\delta_{\nu}=13^2 ST^{\nu}(\mathbf{G}_0)
 =-13 \mathbf{G}_0+\zeta^{\nu} \mathbf{G}_1+\zeta^{2 \nu} \mathbf{G}_2+
  \cdots+\zeta^{12 \nu} \mathbf{G}_{12},\eqno{(3.13)}$$
where the senary sextic forms (i.e., sextic forms in six
variables) are given as follows:
$$\left\{\aligned
  \mathbf{G}_0 =&\mathbf{D}_0^2+\mathbf{D}_{\infty}^2,\\
  \mathbf{G}_1 =&-\mathbf{D}_7^2+2 \mathbf{D}_0 \mathbf{D}_1+10 \mathbf{D}_{\infty}
                 \mathbf{D}_1+2 \mathbf{D}_2 \mathbf{D}_{12}+\\
                &-2 \mathbf{D}_3 \mathbf{D}_{11}-4 \mathbf{D}_4 \mathbf{D}_{10}-2
                 \mathbf{D}_9 \mathbf{D}_5,\\
  \mathbf{G}_2 =&-2 \mathbf{D}_1^2-4 \mathbf{D}_0 \mathbf{D}_2+6 \mathbf{D}_{\infty}
                 \mathbf{D}_2-2 \mathbf{D}_4 \mathbf{D}_{11}+\\
                &+2 \mathbf{D}_5 \mathbf{D}_{10}-2 \mathbf{D}_6 \mathbf{D}_9-2
                 \mathbf{D}_7 \mathbf{D}_8,\\
  \mathbf{G}_3 =&-\mathbf{D}_8^2+2 \mathbf{D}_0 \mathbf{D}_3+10 \mathbf{D}_{\infty}
                 \mathbf{D}_3+2 \mathbf{D}_6 \mathbf{D}_{10}+\\
                &-2 \mathbf{D}_9 \mathbf{D}_7-4 \mathbf{D}_{12} \mathbf{D}_4-2
                 \mathbf{D}_1 \mathbf{D}_2,\\
  \mathbf{G}_4 =&-\mathbf{D}_2^2+10 \mathbf{D}_0 \mathbf{D}_4-2 \mathbf{D}_{\infty}
                 \mathbf{D}_4+2 \mathbf{D}_5 \mathbf{D}_{12}+\\
                &-2 \mathbf{D}_9 \mathbf{D}_8-4 \mathbf{D}_1 \mathbf{D}_3-2
                 \mathbf{D}_{10} \mathbf{D}_7,\\
  \mathbf{G}_5 =&-2 \mathbf{D}_9^2-4 \mathbf{D}_0 \mathbf{D}_5+6 \mathbf{D}_{\infty}
                 \mathbf{D}_5-2 \mathbf{D}_{10} \mathbf{D}_8+\\
                &+2 \mathbf{D}_6 \mathbf{D}_{12}-2 \mathbf{D}_2 \mathbf{D}_3-2
                 \mathbf{D}_{11} \mathbf{D}_7,\\
  \mathbf{G}_6 =&-2 \mathbf{D}_3^2-4 \mathbf{D}_0 \mathbf{D}_6+6 \mathbf{D}_{\infty}
                 \mathbf{D}_6-2 \mathbf{D}_{12} \mathbf{}_7+\\
                &+2 \mathbf{D}_2 \mathbf{D}_4-2 \mathbf{D}_5 \mathbf{D}_1-2
                 \mathbf{D}_8 \mathbf{D}_{11},\\
  \mathbf{G}_7 =&-2 \mathbf{D}_{10}^2+6 \mathbf{D}_0 \mathbf{D}_7+4 \mathbf{D}_{\infty}
                 \mathbf{D}_7-2 \mathbf{D}_1 \mathbf{D}_6+\\
                &-2 \mathbf{D}_2 \mathbf{D}_5-2 \mathbf{D}_8 \mathbf{D}_{12}-2
                 \mathbf{D}_9 \mathbf{D}_{11},\\
  \mathbf{G}_8 =&-2 \mathbf{D}_4^2+6 \mathbf{D}_0 \mathbf{D}_8+4 \mathbf{D}_{\infty}
                 \mathbf{D}_8-2 \mathbf{D}_3 \mathbf{D}_5+\\
                &-2 \mathbf{D}_6 \mathbf{D}_2-2 \mathbf{D}_{11} \mathbf{D}_{10}-2
                 \mathbf{D}_1 \mathbf{D}_7,\\
  \mathbf{G}_9 =&-\mathbf{D}_{11}^2+2 \mathbf{D}_0 \mathbf{D}_9+10 \mathbf{D}_{\infty}
                 \mathbf{D}_9+2 \mathbf{D}_5 \mathbf{D}_4+\\
                &-2 \mathbf{D}_1 \mathbf{D}_8-4 \mathbf{D}_{10} \mathbf{D}_{12}-2
                 \mathbf{D}_3 \mathbf{D}_6,\\
  \mathbf{G}_{10} =&-\mathbf{D}_5^2+10 \mathbf{D}_0 \mathbf{D}_{10}-2 \mathbf{D}_{\infty}
                    \mathbf{D}_{10}+2 \mathbf{D}_6 \mathbf{D}_4+\\
                   &-2 \mathbf{D}_3 \mathbf{D}_7-4 \mathbf{D}_9 \mathbf{D}_1-2
                    \mathbf{D}_{12} \mathbf{D}_{11},\\
  \mathbf{G}_{11} =&-2 \mathbf{D}_{12}^2+6 \mathbf{D}_0 \mathbf{D}_{11}+4 \mathbf{D}_{\infty}
                    \mathbf{D}_{11}-2 \mathbf{D}_9 \mathbf{D}_2+\\
                   &-2 \mathbf{D}_5 \mathbf{D}_6-2 \mathbf{D}_7 \mathbf{D}_4-2
                    \mathbf{D}_3 \mathbf{D}_8,\\
  \mathbf{G}_{12} =&-\mathbf{D}_6^2+10 \mathbf{D}_0 \mathbf{D}_{12}-2 \mathbf{D}_{\infty}
                    \mathbf{D}_{12}+2 \mathbf{D}_2 \mathbf{D}_{10}+\\
                   &-2 \mathbf{D}_1 \mathbf{D}_{11}-4 \mathbf{D}_3 \mathbf{D}_9-2
                    \mathbf{D}_4 \mathbf{D}_8.
\endaligned\right.\eqno{(3.14)}$$
We have that $\mathbf{G}_0$ is invariant under the action of
$\langle H, T \rangle$, a maximal subgroup of order $78$ of
$\text{PSL}(2, 13)$ with index $14$. Note that $\delta_{\infty}$,
$\delta_{\nu}$ for $\nu=0, \ldots, 12$ form an algebraic equation
of degree fourteen. However, we have
$$\Phi_6=\Phi_{0, 1}=\delta_{\infty}+\sum_{\nu=0}^{12} \delta_{\nu}=0.$$
Hence, it is not the Jacobian equation of degree fourteen.

  Recall that the theta functions with characteristic
$\begin{bmatrix} \epsilon\\ \epsilon^{\prime} \end{bmatrix} \in
\mathbb{R}^2$ is defined by the following series which converges
uniformly and absolutely on compact subsets of $\mathbb{C} \times
\mathbb{H}$ (see \cite{FK}, p. 73):
$$\theta \begin{bmatrix} \epsilon\\ \epsilon^{\prime} \end{bmatrix}(z, \tau)
 =\sum_{n \in \mathbb{Z}} \exp \left\{2 \pi i \left[\frac{1}{2}
  \left(n+\frac{\epsilon}{2}\right)^2 \tau+\left(n+\frac{\epsilon}{2}\right)
  \left(z+\frac{\epsilon^{\prime}}{2}\right)\right]\right\}.$$
The modified theta constants (see \cite{FK}, p. 215)
$\varphi_l(\tau):=\theta [\chi_l](0, k \tau)$,
where the characteristic $\chi_l=\begin{bmatrix} \frac{2l+1}{k}\\ 1
\end{bmatrix}$, $l=0, \ldots, \frac{k-3}{2}$, for odd $k$ and
$\chi_l=\begin{bmatrix} \frac{2l}{k}\\ 0 \end{bmatrix}$, $l=0,
\ldots, \frac{k}{2}$, for even $k$. We have the following:

\textbf{Proposition 3.1.} (see \cite{FK}, p. 236). {\it For each
odd integer $k \geq 5$, the map
$\Phi: \tau \mapsto (\varphi_0(\tau), \varphi_1(\tau), \ldots,
 \varphi_{\frac{k-5}{2}}(\tau), \varphi_{\frac{k-3}{2}}(\tau))$
from $\mathbb{H} \cup \mathbb{Q} \cup \{ \infty \}$ to
$\mathbb{C}^{\frac{k-1}{2}}$, defines a holomorphic mapping from
$\overline{\mathbb{H}/\Gamma(k)}$ into
$\mathbb{C} \mathbb{P}^{\frac{k-3}{2}}$.}

  In our case, the map
$\Phi: \tau \mapsto (\varphi_0(\tau), \varphi_1(\tau), \varphi_2(\tau),
 \varphi_3(\tau), \varphi_4(\tau), \varphi_5(\tau))$
gives a holomorphic mapping from the modular curve
$X(13)=\overline{\mathbb{H}/\Gamma(13)}$ into $\mathbb{C} \mathbb{P}^5$,
which corresponds to our six-dimensional representation, i.e., up to
the constants, $z_1, \ldots, z_6$ are just modular forms
$\varphi_0(\tau), \ldots, \varphi_5(\tau)$. Let
$$\left\{\aligned
  a_1(z) &:=e^{-\frac{11 \pi i}{26}} \theta
            \begin{bmatrix} \frac{11}{13}\\ 1 \end{bmatrix}(0, 13z)
           =q^{\frac{121}{104}} \sum_{n \in \mathbb{Z}} (-1)^n
            q^{\frac{1}{2}(13n^2+11n)},\\
  a_2(z) &:=e^{-\frac{7 \pi i}{26}} \theta
            \begin{bmatrix} \frac{7}{13}\\ 1 \end{bmatrix}(0, 13z)
           =q^{\frac{49}{104}} \sum_{n \in \mathbb{Z}} (-1)^n
            q^{\frac{1}{2}(13n^2+7n)},\\
  a_3(z) &:=e^{-\frac{5 \pi i}{26}} \theta
            \begin{bmatrix} \frac{5}{13}\\ 1 \end{bmatrix}(0, 13z)
           =q^{\frac{25}{104}} \sum_{n \in \mathbb{Z}} (-1)^n
            q^{\frac{1}{2}(13n^2+5n)},\\
  a_4(z) &:=-e^{-\frac{3 \pi i}{26}} \theta
            \begin{bmatrix} \frac{3}{13}\\ 1 \end{bmatrix}(0, 13z)
           =-q^{\frac{9}{104}} \sum_{n \in \mathbb{Z}} (-1)^n
            q^{\frac{1}{2}(13n^2+3n)},\\
  a_5(z) &:=e^{-\frac{9 \pi i}{26}} \theta
            \begin{bmatrix} \frac{9}{13}\\ 1 \end{bmatrix}(0, 13z)
           =q^{\frac{81}{104}} \sum_{n \in \mathbb{Z}} (-1)^n
            q^{\frac{1}{2}(13n^2+9n)},\\
  a_6(z) &:=e^{-\frac{\pi i}{26}} \theta
            \begin{bmatrix} \frac{1}{13}\\ 1 \end{bmatrix}(0, 13z)
           =q^{\frac{1}{104}} \sum_{n \in \mathbb{Z}} (-1)^n
            q^{\frac{1}{2}(13n^2+n)}
\endaligned\right.\eqno{(3.15)}$$
be the theta constants of order $13$ and
$$\mathbf{A}(z):=(a_1(z), a_2(z), a_3(z), a_4(z), a_5(z), a_6(z))^{T}.$$
The significance of our six dimensional representation of
$\text{SL}(2, 13)$ comes from the following:

\textbf{Proposition 3.2} (see \cite{Y2}, Proposition 2.5). {\it
If $z \in \mathbb{H}$, then the following relations hold:
$$\mathbf{A}(z+1)=e^{-\frac{3 \pi i}{4}} T \mathbf{A}(z), \quad
  \mathbf{A}\left(-\frac{1}{z}\right)=e^{\frac{\pi i}{4}} \sqrt{z}
  S \mathbf{A}(z),\eqno{(3.16)}$$
where $S$ and $T$ are given in (3.1) and (3.2), and
$0<\text{arg} \sqrt{z} \leq \pi/2$.}

  Note that it is (3.16) that gives an explicit realization of the
isomorphism between the unique sub-representation of parabolic
modular forms of weight $\frac{1}{2}$ on $X(13)$ and the six-dimensional
complex representation of $\text{SL}(2, 13)$ generated by $S$ and $T$.

  Recall that the principal congruence subgroup of level $13$ is the
normal subgroup $\Gamma(13)$ of $\Gamma=\text{SL}(2, \mathbb{Z})$
defined by the exact sequence
$1 \rightarrow \Gamma(13) \rightarrow \Gamma(1) \stackrel{f}{\rightarrow}
 \text{SL}(2, 13) \rightarrow 1$
where $f(\gamma) \equiv \gamma$ (mod $13$) for $\gamma \in \Gamma=\Gamma(1)$.
There is a representation $\rho: \Gamma \rightarrow \text{GL}(6, \mathbb{C})$
with kernel $\Gamma(13)$ defined as follows: if
$t=\begin{pmatrix} 1 & 1\\ 0 & 1 \end{pmatrix}$ and
$s=\begin{pmatrix} 0 & -1\\ 1 & 0 \end{pmatrix}$,
then $\rho(t)=T$ and $\rho(s)=S$. To see that such a representation exists,
note that $\text{PSL}(2, \mathbb{Z}):=\Gamma/\{ \pm I \}$ is defined by
the presentation $\langle s, t; s^2=(st)^3=1 \rangle$ satisfied by $s$ and
$t$ and we have proved that $S$ and $T$ satisfy these relations (in the
projective coordinates). Moreover, we have proved that $\text{PSL}(2, 13)$
is defined by the presentation $\langle S, T; S^2=T^{13}=(ST)^3=1 \rangle$.
Let $p=s t^{-1} s$ and $q=s t^3$. Then
$$h:=q^5 p^2 \cdot p^2 q^6 p^8 \cdot q^5 p^2 \cdot p^3 q
    =\begin{pmatrix}
     4, 428, 249 & -10, 547, 030\\
    -11, 594, 791 & 27, 616, 019
     \end{pmatrix}$$
satisfies that $\rho(h)=H$. The off-diagonal elements of the
matrix $h$, which corresponds to $H$, are congruent to $0$
mod $13$. The connection to $\Gamma_0(13)$ should be obvious.

  Put $x_i(z)=\eta(z) a_i(z)$, $y_i(z)=\eta^3(z) a_i(z)$ and
$u_i(z)=\eta^9(z) a_i(z)$ $(1 \leq i \leq 6)$. Let
$$\left\{\aligned
  X(z) &=(x_1(z), \ldots, x_6(z))^{T},\\
  Y(z) &=(y_1(z), \ldots, y_6(z))^{T}.\\
  U(z) &=(u_1(z), \ldots, u_6(z))^{T}.
\endaligned\right.$$
Then
$$\left\{\aligned
  X(z) &=\eta(z) \mathbf{A}(z)\\
  Y(z) &=\eta^3(z) \mathbf{A}(z),\\
  U(z) &=\eta^9(z) \mathbf{A}(z).
\endaligned\right.$$
Recall that $\eta(z)$ satisfies the following transformation formulas
$\eta(z+1)=e^{\frac{\pi i}{12}} \eta(z)$ and
$\eta\left(-\frac{1}{z}\right)=e^{-\frac{\pi i}{4}} \sqrt{z} \eta(z)$.
By Proposition 3.2, we have
$$X(z+1)=e^{-\frac{2 \pi i}{3}} \rho(t) X(z), \quad
  X\left(-\frac{1}{z}\right)=z \rho(s) X(z),$$
$$Y(z+1)=e^{-\frac{\pi i}{2}} \rho(t) Y(z), \quad
  Y\left(-\frac{1}{z}\right)=e^{-\frac{\pi i}{2}} z^2 \rho(s) Y(z).$$
$$U(z+1)=\rho(t) U(z), \quad
  U\left(-\frac{1}{z}\right)=z^5 \rho(s) U(z).$$
Define $j(\gamma, z):=cz+d$ if $z \in \mathbb{H}$ and
$\gamma=\begin{pmatrix} a & b\\ c & d \end{pmatrix} \in \Gamma(1)$.
Hence,
$$\left\{\aligned
  X(\gamma(z)) &=u(\gamma) j(\gamma, z) \rho(\gamma) X(z),\\
  Y(\gamma(z)) &=v(\gamma) j(\gamma, z)^2 \rho(\gamma) Y(z),\\
  U(\gamma(z)) &=j(\gamma, z)^5 \rho(\gamma) U(z)
\endaligned\right.\eqno{(3.17)}$$
for $\gamma \in \Gamma(1)$, where $u(\gamma)=1, \omega$ or $\omega^2$
with $\omega=e^{\frac{2 \pi i}{3}}$ and $v(\gamma)=\pm 1$ or $\pm i$.
Since $\Gamma(13)=\text{ker}$ $\rho$, we have
$X(\gamma(z))=u(\gamma) j(\gamma, z) X(z)$,
$Y(\gamma(z))=v(\gamma) j(\gamma, z)^2 Y(z)$ and
$U(\gamma(z))=j(\gamma, z)^5 U(z)$ for $\gamma \in \Gamma(13)$.
This means that the functions $x_1(z)$, $\ldots$, $x_6(z)$ are
modular forms of weight one for $\Gamma(13)$ with the same
multiplier $u(\gamma)=1, \omega$ or $\omega^2$ and $y_1(z)$,
$\ldots$, $y_6(z)$ are modular forms of weight two for
$\Gamma(13)$ with the same multiplier $v(\gamma)=\pm 1 $
or $\pm i$.

  From now on, we will use the following abbreviation
$$\mathbf{A}_j=\mathbf{A}_j(a_1(z), \ldots, a_6(z)) \quad
  (0 \leq j \leq 6),$$
$$\mathbf{D}_j=\mathbf{D}_j(a_1(z), \ldots, a_6(z)) \quad
  (j=0,1, \ldots, 12, \infty)$$
and
$$\mathbf{G}_j=\mathbf{G}_j(a_1(z), \ldots, a_6(z)) \quad
  (0 \leq j \leq 12).$$
We have
$$\left\{\aligned
  \mathbf{A}_0 &=q^{\frac{1}{4}} (1+O(q)),\\
  \mathbf{A}_1 &=q^{\frac{17}{52}} (2+O(q)),\\
  \mathbf{A}_2 &=q^{\frac{29}{52}} (2+O(q)),\\
  \mathbf{A}_3 &=q^{\frac{49}{52}} (1+O(q)),\\
  \mathbf{A}_4 &=q^{\frac{25}{52}} (-1+O(q)),\\
  \mathbf{A}_5 &=q^{\frac{9}{52}} (-1+O(q)),\\
  \mathbf{A}_6 &=q^{\frac{1}{52}} (-1+O(q)),
\endaligned\right.$$
and
$$\left\{\aligned
  \mathbf{D}_0 &=q^{\frac{15}{8}} (1+O(q)),\\
  \mathbf{D}_{\infty} &=q^{\frac{7}{8}} (-1+O(q)),\\
  \mathbf{D}_1 &=q^{\frac{99}{104}} (2+O(q)),\\
  \mathbf{D}_2 &=q^{\frac{3}{104}} (-1+O(q)),\\
  \mathbf{D}_3 &=q^{\frac{11}{104}} (1+O(q)),\\
  \mathbf{D}_4 &=q^{\frac{19}{104}} (-2+O(q)),\\
  \mathbf{D}_5 &=q^{\frac{27}{104}} (-1+O(q)),
\endaligned\right. \quad \quad
  \left\{\aligned
  \mathbf{D}_6 &=q^{\frac{35}{104}} (-1+O(q)),\\
  \mathbf{D}_7 &=q^{\frac{43}{104}} (1+O(q)),\\
  \mathbf{D}_8 &=q^{\frac{51}{104}} (3+O(q)),\\
  \mathbf{D}_9 &=q^{\frac{59}{104}} (-2+O(q)),\\
  \mathbf{D}_{10} &=q^{\frac{67}{104}} (1+O(q)),\\
  \mathbf{D}_{11} &=q^{\frac{75}{104}} (-4+O(q)),\\
  \mathbf{D}_{12} &=q^{\frac{83}{104}} (-1+O(q)).
\endaligned\right.$$
Hence,
$$\left\{\aligned
  \mathbf{G}_0 &=q^{\frac{7}{4}} (1+O(q)),\\
  \mathbf{G}_1 &=q^{\frac{43}{52}} (13+O(q)),\\
  \mathbf{G}_2 &=q^{\frac{47}{52}} (-22+O(q)),\\
  \mathbf{G}_3 &=q^{\frac{51}{52}} (-21+O(q)),\\
  \mathbf{G}_4 &=q^{\frac{3}{52}} (-1+O(q)),\\
  \mathbf{G}_5 &=q^{\frac{7}{52}} (2+O(q)),\\
  \mathbf{G}_6 &=q^{\frac{11}{52}} (2+O(q)),
\endaligned\right. \quad \quad
  \left\{\aligned
  \mathbf{G}_7 &=q^{\frac{15}{52}} (-2+O(q)),\\
  \mathbf{G}_8 &=q^{\frac{19}{52}} (-8+O(q)),\\
  \mathbf{G}_9 &=q^{\frac{23}{52}} (6+O(q)),\\
  \mathbf{G}_{10} &=q^{\frac{27}{52}} (1+O(q)),\\
  \mathbf{G}_{11} &=q^{\frac{31}{52}} (-8+O(q)),\\
  \mathbf{G}_{12} &=q^{\frac{35}{52}} (17+O(q)).
\endaligned\right.$$
Note that
$$\aligned
  w_{\nu} &=(\mathbf{A}_0+\zeta^{\nu} \mathbf{A}_1+\zeta^{4 \nu} \mathbf{A}_2
           +\zeta^{9 \nu} \mathbf{A}_3+\zeta^{3 \nu} \mathbf{A}_4
           +\zeta^{12 \nu} \mathbf{A}_5+\zeta^{10 \nu} \mathbf{A}_6)^2\\
          &=\mathbf{A}_0^2+2 (\mathbf{A}_1 \mathbf{A}_5+\mathbf{A}_2 \mathbf{A}_3
           +\mathbf{A}_4 \mathbf{A}_6)+\\
          &+2 \zeta^{\nu} (\mathbf{A}_0 \mathbf{A}_1+\mathbf{A}_2 \mathbf{A}_6)
           +2 \zeta^{3 \nu} (\mathbf{A}_0 \mathbf{A}_4+\mathbf{A}_2 \mathbf{A}_5)+\\
          &+2 \zeta^{9 \nu} (\mathbf{A}_0 \mathbf{A}_3+\mathbf{A}_5 \mathbf{A}_6)
           +2 \zeta^{12 \nu} (\mathbf{A}_0 \mathbf{A}_5+\mathbf{A}_3 \mathbf{A}_4)+\\
          &+2 \zeta^{10 \nu} (\mathbf{A}_0 \mathbf{A}_6+\mathbf{A}_1 \mathbf{A}_3)
           +2 \zeta^{4 \nu} (\mathbf{A}_0 \mathbf{A}_2+\mathbf{A}_1 \mathbf{A}_4)+\\
          &+\zeta^{2 \nu} (\mathbf{A}_1^2+2 \mathbf{A}_4 \mathbf{A}_5)
           +\zeta^{5 \nu} (\mathbf{A}_3^2+2 \mathbf{A}_1 \mathbf{A}_2)+\\
          &+\zeta^{6 \nu} (\mathbf{A}_4^2+2 \mathbf{A}_3 \mathbf{A}_6)
           +\zeta^{11 \nu} (\mathbf{A}_5^2+2 \mathbf{A}_1 \mathbf{A}_6)+\\
          &+\zeta^{8 \nu} (\mathbf{A}_2^2+2 \mathbf{A}_3 \mathbf{A}_5)
           +\zeta^{7 \nu} (\mathbf{A}_6^2+2 \mathbf{A}_4 \mathbf{A}_2),
\endaligned$$
where
$$\mathbf{A}_0^2+2 (\mathbf{A}_1 \mathbf{A}_5+\mathbf{A}_2 \mathbf{A}_3
  +\mathbf{A}_4 \mathbf{A}_6)=q^{\frac{1}{2}} (-1+O(q)),$$
$$\left\{\aligned
  \mathbf{A}_0 \mathbf{A}_1+\mathbf{A}_2 \mathbf{A}_6 &=q^{\frac{41}{26}} (-3+O(q)),\\
  \mathbf{A}_0 \mathbf{A}_4+\mathbf{A}_2 \mathbf{A}_5 &=q^{\frac{19}{26}} (-3+O(q)),\\
  \mathbf{A}_0 \mathbf{A}_3+\mathbf{A}_5 \mathbf{A}_6 &=q^{\frac{5}{26}} (1+O(q)),\\
  \mathbf{A}_0 \mathbf{A}_5+\mathbf{A}_3 \mathbf{A}_4 &=q^{\frac{11}{26}} (-1+O(q)),\\
  \mathbf{A}_0 \mathbf{A}_6+\mathbf{A}_1 \mathbf{A}_3 &=q^{\frac{7}{26}} (-1+O(q)),\\
  \mathbf{A}_0 \mathbf{A}_2+\mathbf{A}_1 \mathbf{A}_4 &=q^{\frac{47}{26}} (-1+O(q)),
\endaligned\right.$$
and
$$\left\{\aligned
  \mathbf{A}_1^2+2 \mathbf{A}_4 \mathbf{A}_5 &=q^{\frac{17}{26}} (6+O(q)),\\
  \mathbf{A}_3^2+2 \mathbf{A}_1 \mathbf{A}_2 &=q^{\frac{23}{26}} (8+O(q)),\\
  \mathbf{A}_4^2+2 \mathbf{A}_3 \mathbf{A}_6 &=q^{\frac{25}{26}} (-1+O(q)),\\
  \mathbf{A}_5^2+2 \mathbf{A}_1 \mathbf{A}_6 &=q^{\frac{9}{26}} (-3+O(q)),\\
  \mathbf{A}_2^2+2 \mathbf{A}_3 \mathbf{A}_5 &=q^{\frac{29}{26}} (2+O(q)),\\
  \mathbf{A}_6^2+2 \mathbf{A}_4 \mathbf{A}_2 &=q^{\frac{1}{26}} (1+O(q)).
\endaligned\right.$$

  In the introduction, the invariant homogeneous polynomials $\Phi_{m, n}$
are defined. Now we give the normalization for the following four families
of polynomials of degrees $d=12$, $16$, $20$ and $30$. For $d=12$, there are
two such invariant homogeneous polynomials $\Phi_{3, 0}$ and $\Phi_{0, 2}$:
$$\left\{\aligned
  \Phi_{3, 0} &:=-\frac{1}{13 \cdot 30} \left(\sum_{\nu=0}^{12}
                   w_{\nu}^3+w_{\infty}^3\right),\\
  \Phi_{0, 2} &:=-\frac{1}{13 \cdot 52} \left(\sum_{\nu=0}^{12}
                 \delta_{\nu}^2+\delta_{\infty}^2\right),
\endaligned\right.\eqno{(3.18)}$$
For $d=16$, there are two such invariant homogeneous polynomials
$\Phi_{4, 0}$ and $\Phi_{1, 2}$ which need not to be normalized.
For $d=20$, there are two such invariant homogeneous polynomials
$\Phi_{5, 0}$ and $\Phi_{2, 2}$:
$$\left\{\aligned
  \Phi_{5, 0} &:=\frac{1}{13 \cdot 25} \left(\sum_{\nu=0}^{12} w_{\nu}^5
                 +w_{\infty}^5\right),\\
  \Phi_{2, 2} &:=\frac{1}{13 \cdot 26} \left(\sum_{\nu=0}^{12} w_{\nu}^2
                 \delta_{\nu}^2+w_{\infty}^2 \delta_{\infty}^2\right)
\endaligned\right.\eqno{(3.19)}$$
For $d=30$, there are three such invariant homogeneous polynomials
$\Phi_{0, 5}$, $\Phi_{3, 3}$ and $\Phi_{6, 1}$:
$$\left\{\aligned
  \Phi_{0, 5} &:=-\frac{1}{13 \cdot 1315} \left(\sum_{\nu=0}^{12}
                 \delta_{\nu}^5+\delta_{\infty}^5\right),\\
  \Phi_{3, 3} &:=-\frac{1}{13 \cdot 27} \left(\sum_{\nu=0}^{12}
                 w_{\nu}^3 \delta_{\nu}^3+w_{\infty}^3 \delta_{\infty}^3\right),\\
  \Phi_{6, 1} &:=-\frac{1}{13 \cdot 285} \left(\sum_{\nu=0}^{12}
                 w_{\nu}^6 \delta_{\nu}+w_{\infty}^6 \delta_{\infty}\right),
\endaligned\right.\eqno{(3.20)}$$

\textbf{Theorem 3.3.} {\it The $G$-invariant homogeneous polynomials
$\Phi_{m, n}$ of degrees $d=4$, $8$, $10$, $12$, $14$, $16$, $20$ and
$30$ in $x_1(z)$, $\ldots$, $x_6(z)$ can be identified with modular
forms as follows$:$
$$\left\{\aligned
  \Phi_{4}(x_1(z), \ldots, x_6(z)) &=0,\\
  \Phi_{8}(x_1(z), \ldots, x_6(z)) &=0,\\
  \Phi_{10}(x_1(z), \ldots, x_6(z)) &=0,\\
  \Phi_{3, 0}(x_1(z), \ldots, x_6(z)) &=\Delta(z),\\
  \Phi_{0, 2}(x_1(z), \ldots, x_6(z)) &=\Delta(z),\\
  \Phi_{14}(x_1(z), \ldots, x_6(z)) &=0,\\
  \Phi_{4, 0}(x_1(z), \ldots, x_6(z)) &=0,\\
  \Phi_{1, 2}(x_1(z), \ldots, x_6(z)) &=0,\\
  \Phi_{5, 0}(x_1(z), \ldots, x_6(z)) &=\eta(z)^8 \Delta(z) E_4(z),\\
  \Phi_{2, 2}(x_1(z), \ldots, x_6(z)) &=\eta(z)^8 \Delta(z) E_4(z),\\
  \Phi_{0, 5}(x_1(z), \ldots, x_6(z)) &=\Delta(z)^2 E_6(z),\\
  \Phi_{3, 3}(x_1(z), \ldots, x_6(z)) &=\Delta(z)^2 E_6(z),\\
  \Phi_{6, 1}(x_1(z), \ldots, x_6(z)) &=\Delta(z)^2 E_6(z).\\
\endaligned\right.\eqno{(3.21)}$$}

{\it Proof}. We divide the proof into four parts (see also \cite{Y3}).
The first part is the calculation of $\Phi_{5, 0}$ and $\Phi_{3, 0}$.
Up to a constant, $\Phi_{5, 0}=w_0^5+w_1^5+\cdots+w_{12}^5+w_{\infty}^5$.
As a polynomial in six variables, $\Phi_{5, 0}(z_1, z_2, z_3, z_4, z_5, z_6)$
is a $G$-invariant polynomial. Moreover, for $\gamma \in \Gamma(1)$,
$$\aligned
  &\Phi_{5, 0}(Y(\gamma(z))^{T})=\Phi_{5, 0}(v(\gamma)
   j(\gamma, z)^2 (\rho(\gamma) Y(z))^{T})\\
 =&v(\gamma)^{20} j(\gamma, z)^{40} \Phi_{5, 0}((\rho(\gamma) Y(z))^{T})
 =j(\gamma, z)^{40} \Phi_{5, 0}((\rho(\gamma) Y(z))^{T}).
\endaligned$$
Note that $\rho(\gamma) \in \langle \rho(s), \rho(t) \rangle=G$ and
$\Phi_{5, 0}$ is a $G$-invariant polynomial, we have
$$\Phi_{5, 0}(Y(\gamma(z))^{T})=j(\gamma, z)^{40} \Phi_{5, 0}(Y(z)^{T}),
  \quad \text{for $\gamma \in \Gamma(1)$}.$$
This implies that $\Phi_{5, 0}(y_1(z), \ldots, y_6(z))$ is a modular form
of weight $40$ for the full modular group $\Gamma(1)$. Moreover, we will
show that it is a cusp form. In fact,
$$\aligned
  &\Phi_{5, 0}(a_1(z), \ldots, a_6(z))=13^5 q^{\frac{5}{2}} (1+O(q))^5+\\
  &+\sum_{\nu=0}^{12} [q^{\frac{1}{2}} (-1+O(q))+\\
  &+2 \zeta^{\nu} q^{\frac{41}{26}} (-3+O(q))+2 \zeta^{3 \nu}
    q^{\frac{19}{26}} (-3+O(q))
   +2 \zeta^{9 \nu} q^{\frac{5}{26}} (1+O(q))+\\
  &+2 \zeta^{12 \nu} q^{\frac{11}{26}} (-1+O(q))+2 \zeta^{10 \nu}
   q^{\frac{7}{26}} (-1+O(q))
   +2 \zeta^{4 \nu} q^{\frac{47}{26}} (-1+O(q))+\\
  &+\zeta^{2 \nu} q^{\frac{17}{26}} (6+O(q))+\zeta^{5 \nu}
   q^{\frac{23}{26}} (8+O(q))
   +\zeta^{6 \nu} q^{\frac{25}{26}} (-1+O(q))+\\
  &+\zeta^{11 \nu} q^{\frac{9}{26}} (-3+O(q))+\zeta^{8 \nu}
   q^{\frac{29}{26}} (2+O(q))
   +\zeta^{7 \nu} q^{\frac{1}{26}} (1+O(q))]^5.
\endaligned$$
We will calculate the $q^{\frac{1}{2}}$-term which is the lowest
degree. For the partition $13=4 \cdot 1+9$, the corresponding
term is
$$\begin{pmatrix} 5\\ 4, 1 \end{pmatrix} (\zeta^{7 \nu} q^{\frac{1}{26}})^4
  \cdot (-3) \zeta^{11 \nu} q^{\frac{9}{26}}=-15 q^{\frac{1}{2}}.$$
For the partition $13=3 \cdot 1+2 \cdot 5$, the corresponding term
is
$$\begin{pmatrix} 5\\ 3, 2 \end{pmatrix} (\zeta^{7 \nu} q^{\frac{1}{26}})^3
  \cdot (2 \zeta^{9 \nu} q^{\frac{5}{26}})^2=40 q^{\frac{1}{2}}.$$
Hence, for $\Phi_{5, 0}(y_1(z), \ldots, y_6(z))$ which is a modular
form for $\Gamma(1)$ with weight $40$, the lowest degree term is given by
$$(-15+40) q^{\frac{1}{2}} \cdot q^{\frac{3}{24} \cdot 20}=25 q^3.$$
Thus,
$$\Phi_{5, 0}(y_1(z), \ldots, y_6(z))=q^3 (13 \cdot 25+O(q)).$$
The leading term of $\Phi_{5, 0}(y_1(z), \ldots, y_6(z))$ together
with its weight $40$ suffice to identify this modular form with
$\Phi_{5, 0}(y_1(z), \ldots, y_6(z))=13 \cdot 25 \Delta(z)^3 E_4(z)$.
Consequently,
$$\Phi_{5, 0}(x_1(z), \ldots, x_6(z))=13 \cdot 25 \Delta(z)^3
  E_4(z)/\eta(z)^{40}=13 \cdot 25 \eta(z)^8 \Delta(z) E_4(z).$$

  Up to a constant, $\Phi_{3, 0}=w_0^3+w_1^3+\cdots+w_{12}^3+w_{\infty}^3$,
The calculation of $\Phi_{3, 0}$ is similar as that of $\Phi_{5, 0}$.
We find that
$$\Phi_{3, 0}(x_1(z), \ldots, x_6(z))=-13 \cdot 30 \Delta(z).$$

  The second part is the calculation of $\Phi_4$, $\Phi_8$ and
$\Phi_{4, 0}$. The calculation of $\Phi_4$ has been done in \cite{Y2},
Theorem 3.1. We will give the calculation of $\Phi_{4, 0}$. Up to a
constant, $\Phi_{4, 0}=w_0^4+w_1^4+\cdots+w_{12}^4+w_{\infty}^4$.
Similar as the above calculation for $\Phi_{5, 0}$, we find that
$\Phi_{4, 0}(y_1(z), \ldots, y_6(z))$ is a modular form of weight
$32$ for the full modular group $\Gamma(1)$. Moreover, we will
show that it is a cusp form. In fact,
$$\aligned
  &\Phi_{4, 0}(a_1(z), \ldots, a_6(z))=13^4 q^2 (1+O(q))^4+\\
  &+\sum_{\nu=0}^{12} [q^{\frac{1}{2}} (-1+O(q))+\\
  &+2 \zeta^{\nu} q^{\frac{41}{26}} (-3+O(q))+2 \zeta^{3 \nu} q^{\frac{19}{26}} (-3+O(q))
   +2 \zeta^{9 \nu} q^{\frac{5}{26}} (1+O(q))+\\
  &+2 \zeta^{12 \nu} q^{\frac{11}{26}} (-1+O(q))+2 \zeta^{10 \nu} q^{\frac{7}{26}} (-1+O(q))
   +2 \zeta^{4 \nu} q^{\frac{47}{26}} (-1+O(q))+\\
  &+\zeta^{2 \nu} q^{\frac{17}{26}} (6+O(q))+\zeta^{5 \nu} q^{\frac{23}{26}} (8+O(q))
   +\zeta^{6 \nu} q^{\frac{25}{26}} (-1+O(q))+\\
  &+\zeta^{11 \nu} q^{\frac{9}{26}} (-3+O(q))+\zeta^{8 \nu} q^{\frac{29}{26}} (2+O(q))
   +\zeta^{7 \nu} q^{\frac{1}{26}} (1+O(q))]^4.
\endaligned$$
We will calculate the $q$-term which is the lowest degree. For example,
consider the partition $26=3 \cdot 1+23$, the corresponding term is
$$\begin{pmatrix} 4\\ 3, 1 \end{pmatrix} (\zeta^{7 \nu} q^{\frac{1}{26}})^3
  \cdot 8 \zeta^{5 \nu} q^{\frac{23}{26}}=32 q.$$
For the other partitions, the calculation is similar. In conclusion,
we find that the coefficients of the $q$-term is an integer. Hence,
for $\Phi_{4, 0}(y_1(z), \ldots, y_6(z))$ which is a modular form for
$\Gamma(1)$ with weight $32$, the lowest degree term is given by
$$\text{some integer} \cdot q \cdot q^{\frac{3}{24} \cdot 16}
 =\text{some integer} \cdot q^3.$$
This implies that $\Phi_{4, 0}(y_1(z), \ldots, y_6(z))$ has a factor of
$\Delta(z)^3$, which is a cusp form of weight $36$. Therefore,
$\Phi_{4, 0}(y_1(z), \ldots, y_6(z))=0$. The calculation of $\Phi_8$ is
similar as that of $\Phi_{4, 0}$.

   The third part is the calculation of $\Phi_{0, 2}$ and $\Phi_{0, 5}$.
Up to a constant,
$\Phi_{0, 2}=\delta_0^2+\delta_1^2+\cdots+\delta_{12}^2+\delta_{\infty}^2$.
As a polynomial in six variables, $\Phi_{0, 2}(z_1, z_2, z_3, z_4, z_5, z_6)$
is a $G$-invariant polynomial. Moreover, for $\gamma \in \Gamma(1)$,
$$\aligned
  &\Phi_{0, 2}(X(\gamma(z))^{T})
 =\Phi_{0, 2}(u(\gamma) j(\gamma, z) (\rho(\gamma) X(z))^{T})\\
 =&u(\gamma)^{12} j(\gamma, z)^{12} \Phi_{0, 2}((\rho(\gamma) X(z))^{T})
 =j(\gamma, z)^{12} \Phi_{0, 2}((\rho(\gamma) X(z))^{T}).
\endaligned$$
Note that $\rho(\gamma) \in \langle \rho(s), \rho(t) \rangle=G$
and $\Phi_{0, 2}$ is a $G$-invariant polynomial, we have
$$\Phi_{0, 2}(X(\gamma(z))^{T})=j(\gamma, z)^{12}
  \Phi_{0, 2}(X(z)^{T}), \quad \text{for $\gamma \in
  \Gamma(1)$}.$$
This implies that $\Phi_{0, 2}(x_1(z), \ldots, x_6(z))$ is a modular
form of weight $12$ for the full modular group $\Gamma(1)$. Moreover,
we will show that it is a cusp form. In fact,
$$\aligned
  &\Phi_{0, 2}(a_1(z), \ldots, a_6(z))=13^4 q^{\frac{7}{2}} (1+O(q))^2+\\
  &+\sum_{\nu=0}^{12} [-13 q^{\frac{7}{4}} (1+O(q))+\\
  &+\zeta^{\nu} q^{\frac{43}{52}} (13+O(q))+\zeta^{2 \nu}
   q^{\frac{47}{52}} (-22+O(q))
   +\zeta^{3 \nu} q^{\frac{51}{52}} (-21+O(q))+\\
  &+\zeta^{4 \nu} q^{\frac{3}{52}} (-1+O(q))+\zeta^{5 \nu}
   q^{\frac{7}{52}} (2+O(q))
   +\zeta^{6 \nu} q^{\frac{11}{52}} (2+O(q))+\\
  &+\zeta^{7 \nu} q^{\frac{15}{52}} (-2+O(q))+\zeta^{8 \nu}
   q^{\frac{19}{52}} (-8+O(q))
   +\zeta^{9 \nu} q^{\frac{23}{52}} (6+O(q))+\\
  &+\zeta^{10 \nu} q^{\frac{27}{52}} (1+O(q))+\zeta^{11 \nu}
   q^{\frac{31}{52}} (-8+O(q))
   +\zeta^{12 \nu} q^{\frac{35}{52}} (17+O(q))]^2.
\endaligned$$
We will calculate the $q^{\frac{1}{2}}$-term which is the lowest
degree. For the partition $26=3+23$, the corresponding term is
$$\begin{pmatrix} 2\\ 1, 1 \end{pmatrix} \zeta^{4 \nu} q^{\frac{3}{52}}
  \cdot (-1) \cdot \zeta^{9 \nu} q^{\frac{23}{52}} \cdot 6=-12 q^{\frac{1}{2}}.$$
For the partition $26=7+19$, the corresponding term is
$$\begin{pmatrix} 2\\ 1, 1 \end{pmatrix} \zeta^{5 \nu} q^{\frac{7}{52}}
  \cdot 2 \cdot \zeta^{8 \nu} q^{\frac{19}{52}} \cdot (-8)=-32 q^{\frac{1}{2}}.$$
For the partition $26=11+15$, the corresponding term is
$$\begin{pmatrix} 2\\ 1, 1 \end{pmatrix} \zeta^{6 \nu} q^{\frac{11}{52}}
  \cdot 2 \cdot \zeta^{7 \nu} q^{\frac{15}{52}} \cdot (-2)=-8 q^{\frac{1}{2}}.$$
Hence, for $\Phi_{0, 2}(x_1(z), \ldots, x_6(z))$ which is a
modular form for $\Gamma(1)$ with weight $12$, the lowest degree
term is given by $(-12-32-8) q^{\frac{1}{2}} \cdot q^{\frac{12}{24}}=-52 q$.
Thus,
$$\Phi_{0, 2}(x_1(z), \ldots, x_6(z))=q (-13 \cdot 52+O(q)).$$
The leading term of $\Phi_{0, 2}(x_1(z), \ldots, x_6(z))$ together
with its weight $12$ suffice to identify this modular form with
$$\Phi_{0, 2}(x_1(z), \ldots, x_6(z))=-13 \cdot 52 \Delta(z).$$

  Up to a constant,
$\Phi_{0, 5}=\delta_0^5+\delta_1^5+\cdots+\delta_{12}^5+\delta_{\infty}^5$.
As a polynomial in six variables, $\Phi_{0, 5}(z_1, z_2, z_3, z_4, z_5, z_6)$
is a $G$-invariant polynomial. Similarly as above, we can show that
$\Phi_{0, 5}(x_1(z), \ldots, x_6(z))$ is a modular form of weight $30$ for
the full modular group $\Gamma(1)$. Moreover, we will show that it is a
cusp form. In fact,
$$\aligned
  &\Phi_{0, 5}(a_1(z), \ldots, a_6(z))=13^{10} q^{\frac{35}{4}} (1+O(q))^5+\\
  &+\sum_{\nu=0}^{12} [-13 q^{\frac{7}{4}} (1+O(q))+\\
  &+\zeta^{\nu} q^{\frac{43}{52}} (13+O(q))+\zeta^{2 \nu} q^{\frac{47}{52}} (-22+O(q))
   +\zeta^{3 \nu} q^{\frac{51}{52}} (-21+O(q))+\\
  &+\zeta^{4 \nu} q^{\frac{3}{52}} (-1+O(q))+\zeta^{5 \nu} q^{\frac{7}{52}} (2+O(q))
   +\zeta^{6 \nu} q^{\frac{11}{52}} (2+O(q))+\\
  &+\zeta^{7 \nu} q^{\frac{15}{52}} (-2+O(q))+\zeta^{8 \nu} q^{\frac{19}{52}} (-8+O(q))
   +\zeta^{9 \nu} q^{\frac{23}{52}} (6+O(q))+\\
  &+\zeta^{10 \nu} q^{\frac{27}{52}} (1+O(q))+\zeta^{11 \nu} q^{\frac{31}{52}} (-8+O(q))
   +\zeta^{12 \nu} q^{\frac{35}{52}} (17+O(q))]^5.
\endaligned$$
We will calculate the $q^{\frac{3}{4}}$-term which is the lowest degree.
(1) For the partition $39=4 \cdot 3+27$, the corresponding term is
$$\begin{pmatrix} 5\\ 4, 1 \end{pmatrix} (\zeta^{4 \nu}
  q^{\frac{3}{52}} \cdot (-1))^4 \cdot \zeta^{10 \nu}
  q^{\frac{27}{52}}=5 q^{\frac{3}{4}}.$$
(2) For the partition $39=3 \cdot 3+7+23$, the corresponding term is
$$\begin{pmatrix} 5\\ 3, 1, 1 \end{pmatrix} (\zeta^{4 \nu}
  q^{\frac{3}{52}} \cdot (-1))^3 \cdot \zeta^{5 \nu} q^{\frac{7}{52}}
  \cdot 2 \cdot \zeta^{9 \nu} q^{\frac{23}{52}} \cdot 6=-240 q^{\frac{3}{4}}.$$
(3) For the partition $39=3 \cdot 3+11+19$, the corresponding term is
$$\begin{pmatrix} 5\\ 3, 1, 1 \end{pmatrix} (\zeta^{4 \nu}
  q^{\frac{3}{52}} \cdot (-1))^3 \cdot
  \zeta^{6 \nu} q^{\frac{11}{52}} \cdot 2 \cdot \zeta^{8 \nu}
  q^{\frac{19}{52}} \cdot (-8)=320 q^{\frac{3}{4}}.$$
(4) For the partition $39=3 \cdot 3+2 \cdot 15$, the corresponding term is
$$\begin{pmatrix} 5\\ 3, 2 \end{pmatrix} (\zeta^{4 \nu}
  q^{\frac{3}{52}} \cdot (-1))^3 \cdot
  (\zeta^{7 \nu} q^{\frac{15}{52}} \cdot (-2))^2=-40 q^{\frac{3}{4}}.$$
(5) For the partition $39=2 \cdot 3+3 \cdot 11$, the corresponding term is
$$\begin{pmatrix} 5\\ 2, 3 \end{pmatrix} (\zeta^{4 \nu}
  q^{\frac{3}{52}} \cdot (-1))^2 \cdot
  (\zeta^{6 \nu} q^{\frac{11}{52}} \cdot 2)^3=80 q^{\frac{3}{4}}.$$
(6) For the partition $39=2 \cdot 3+2 \cdot 7+19$, the corresponding term is
$$\begin{pmatrix} 5\\ 2, 2, 1 \end{pmatrix} (\zeta^{4 \nu}
  q^{\frac{3}{52}} \cdot (-1))^2 \cdot
  (\zeta^{5 \nu} q^{\frac{7}{52}} \cdot 2)^2 \cdot \zeta^{8 \nu}
  q^{\frac{19}{52}} \cdot (-8)=-960 q^{\frac{3}{4}}.$$
(7) For the partition $39=2 \cdot 3+7+11+15$, the corresponding term is
$$\begin{pmatrix} 5\\ 2, 1, 1, 1 \end{pmatrix} (\zeta^{4 \nu}
  q^{\frac{3}{52}} \cdot (-1))^2 \cdot
  \zeta^{5 \nu} q^{\frac{7}{52}} \cdot 2 \cdot \zeta^{6 \nu}
  q^{\frac{11}{52}} \cdot 2 \cdot \zeta^{7 \nu}
  q^{\frac{15}{52}} \cdot (-2)=-480 q^{\frac{3}{4}}.$$
(8) For the partition $39=1 \cdot 3+3 \cdot 7+15$, the corresponding
term is
$$\begin{pmatrix} 5\\ 1, 3, 1 \end{pmatrix} \zeta^{4 \nu}
  q^{\frac{3}{52}} \cdot (-1) \cdot
  (\zeta^{5 \nu} q^{\frac{7}{52}} \cdot 2)^3 \cdot \zeta^{7 \nu}
  q^{\frac{15}{52}} \cdot (-2)=320 q^{\frac{3}{4}}.$$
(9) For the partition $39=1 \cdot 3+2 \cdot 7+2 \cdot 11$, the
corresponding term is
$$\begin{pmatrix} 5\\ 1, 2, 2 \end{pmatrix} \zeta^{4 \nu}
  q^{\frac{3}{52}} \cdot (-1) \cdot
  (\zeta^{5 \nu} q^{\frac{7}{52}} \cdot 2)^2 \cdot (\zeta^{6 \nu}
  q^{\frac{11}{52}} \cdot 2)^2=-480 q^{\frac{3}{4}}.$$
(10) For the partition $39=4 \cdot 7+11$, the corresponding term is
$$\begin{pmatrix} 5\\ 4, 1 \end{pmatrix} (\zeta^{5 \nu}
  q^{\frac{7}{52}} \cdot 2)^4 \cdot
  \zeta^{6 \nu} q^{\frac{11}{52}} \cdot 2=160 q^{\frac{3}{4}}.$$
Hence, for $\Phi_{0, 5}(x_1(z), \ldots, x_6(z))$ which is a modular
form for $\Gamma(1)$ with weight $30$, the lowest degree term is
given by
$$(5-240+320-40+80-960-480+320-480+160) q^{\frac{3}{4}} \cdot
  q^{\frac{30}{24}}=-1315 q^2.$$
Thus,
$$\Phi_{0, 5}(x_1(z), \ldots, x_6(z))=q^2 (-13 \cdot 1315+O(q)).$$
The leading term of $\Phi_{0, 5}(x_1(z), \ldots, x_6(z))$ together
with its weight $30$ suffice to identify this modular form with
$$\Phi_{0, 5}(x_1(z), \ldots, x_6(z))=-13 \cdot 1315 \Delta(z)^2 E_6(z).$$

  The last part is the calculation of $\Phi_{10}$, $\Phi_{14}$, $\Phi_{1, 2}$,
$\Phi_{2, 2}$, $\Phi_{3, 3}$ and $\Phi_{6, 1}$. For
$\Phi_{10}=w_0 \delta_0+w_1 \delta_1+\cdots+w_{12} \delta_{12}+w_{\infty} \delta_{\infty}$.
As a polynomial in six variables, $\Phi_{10}(z_1, z_2, z_3, z_4, z_5, z_6)$
is a $G$-invariant polynomial. Moreover, for $\gamma \in \Gamma(1)$,
$$\Phi_{10}(U(\gamma(z))^{T})=\Phi_{10}(j(\gamma, z)^5 (\rho(\gamma) U(z))^{T})
 =j(\gamma, z)^{50} \Phi_{10}((\rho(\gamma) U(z))^{T}).$$
Note that $\rho(\gamma) \in \langle \rho(s), \rho(t) \rangle=G$ and
$\Phi_{10}$ is a $G$-invariant polynomial, we have
$$\Phi_{10}(U(\gamma(z))^{T})=j(\gamma, z)^{50} \Phi_{10}(U(z)^{T}),
  \quad \text{for $\gamma \in \Gamma(1)$}.$$
This implies that $\Phi_{10}(u_1(z), \ldots, u_6(z))$ is a modular form
of weight $50$ for the full modular group $\Gamma(1)$. Moreover, we will
show that it is a cusp form. In fact,
$$\aligned
  &\Phi_{10}(a_1(z), \ldots, a_6(z))=13^3 q^{\frac{9}{4}} (1+O(q))^2+\\
  &+\sum_{\nu=0}^{12} [q^{\frac{1}{2}} (-1+O(q))+\\
  &+2 \zeta^{\nu} q^{\frac{41}{26}} (-3+O(q))+2 \zeta^{3 \nu}
    q^{\frac{19}{26}} (-3+O(q))
   +2 \zeta^{9 \nu} q^{\frac{5}{26}} (1+O(q))+\\
  &+2 \zeta^{12 \nu} q^{\frac{11}{26}} (-1+O(q))+2 \zeta^{10 \nu}
   q^{\frac{7}{26}} (-1+O(q))
   +2 \zeta^{4 \nu} q^{\frac{47}{26}} (-1+O(q))+\\
  &+\zeta^{2 \nu} q^{\frac{17}{26}} (6+O(q))+\zeta^{5 \nu}
   q^{\frac{23}{26}} (8+O(q))
   +\zeta^{6 \nu} q^{\frac{25}{26}} (-1+O(q))+\\
  &+\zeta^{11 \nu} q^{\frac{9}{26}} (-3+O(q))+\zeta^{8 \nu}
   q^{\frac{29}{26}} (2+O(q))
   +\zeta^{7 \nu} q^{\frac{1}{26}} (1+O(q))]\\
  &\times [-13 q^{\frac{7}{4}} (1+O(q))+\\
  &+\zeta^{\nu} q^{\frac{43}{52}} (13+O(q))+\zeta^{2 \nu}
   q^{\frac{47}{52}} (-22+O(q))
   +\zeta^{3 \nu} q^{\frac{51}{52}} (-21+O(q))+\\
  &+\zeta^{4 \nu} q^{\frac{3}{52}} (-1+O(q))+\zeta^{5 \nu}
   q^{\frac{7}{52}} (2+O(q))
   +\zeta^{6 \nu} q^{\frac{11}{52}} (2+O(q))+\\
  &+\zeta^{7 \nu} q^{\frac{15}{52}} (-2+O(q))+\zeta^{8 \nu}
   q^{\frac{19}{52}} (-8+O(q))
   +\zeta^{9 \nu} q^{\frac{23}{52}} (6+O(q))+\\
  &+\zeta^{10 \nu} q^{\frac{27}{52}} (1+O(q))+\zeta^{11 \nu}
   q^{\frac{31}{52}} (-8+O(q))
   +\zeta^{12 \nu} q^{\frac{35}{52}} (17+O(q))]
\endaligned$$
We will calculate the $q^{\frac{1}{4}}$-term which is the lowest
degree:
$$2 q^{\frac{5}{26}} \cdot q^{\frac{3}{52}} \cdot (-1)+
  q^{\frac{1}{26}} \cdot 1 \cdot q^{\frac{11}{52}} \cdot 2=0.$$
Hence, for $\Phi_{10}(u_1(z), \ldots, u_6(z))$ which is a modular form for
$\Gamma(1)$ with weight $50$, the lowest degree term is given by
$$\text{some integer} \cdot q^{\frac{5}{4}} \cdot q^{\frac{9}{24} \cdot 10}
 =\text{some integer} \cdot q^5.$$
This implies that $\Phi_{10}(u_1(z), \ldots, u_6(z))$ has a factor of
$\Delta(z)^5$, which is a cusp form of weight $60$. Therefore,
$\Phi_{10}(u_1(z), \ldots, u_6(z))=0$. Consequently,
$\Phi_{10}(x_1(z), \ldots, x_6(z))=0$. The calculation of $\Phi_{14}$
and $\Phi_{1, 2}$ is similar as that of $\Phi_{10}$.

  For
$\Phi_{2, 2}=w_0^2 \delta_0^2+w_1^2 \delta_1^2+\cdots+w_{12}^2 \delta_{12}^2+
w_{\infty}^2 \delta_{\infty}^2$.
As a polynomial in six variables, $\Phi_{2, 2}(z_1, z_2, z_3, z_4, z_5, z_6)$
is a $G$-invariant polynomial. Moreover, for $\gamma \in \Gamma(1)$,
$$\Phi_{2, 2}(U(\gamma(z))^{T})=\Phi_{2, 2}(j(\gamma, z)^5 (\rho(\gamma) U(z))^{T})
 =j(\gamma, z)^{100} \Phi_{2, 2}((\rho(\gamma) U(z))^{T}).$$
Note that $\rho(\gamma) \in \langle \rho(s), \rho(t) \rangle=G$ and
$\Phi_{2, 2}$ is a $G$-invariant polynomial, we have
$$\Phi_{2, 2}(U(\gamma(z))^{T})=j(\gamma, z)^{100} \Phi_{2, 2}(U(z)^{T}),
  \quad \text{for $\gamma \in \Gamma(1)$}.$$
This implies that $\Phi_{2, 2}(u_1(z), \ldots, u_6(z))$ is a modular form
of weight $100$ for the full modular group $\Gamma(1)$. Moreover, we will
show that it is a cusp form. In fact,
$$\aligned
  &\Phi_{2, 2}(a_1(z), \ldots, a_6(z))=13^6 q^{\frac{9}{2}} (1+O(q))^2+\\
  &+\sum_{\nu=0}^{12} [q^{\frac{1}{2}} (-1+O(q))+\\
  &+2 \zeta^{\nu} q^{\frac{41}{26}} (-3+O(q))+2 \zeta^{3 \nu}
    q^{\frac{19}{26}} (-3+O(q))
   +2 \zeta^{9 \nu} q^{\frac{5}{26}} (1+O(q))+\\
  &+2 \zeta^{12 \nu} q^{\frac{11}{26}} (-1+O(q))+2 \zeta^{10 \nu}
   q^{\frac{7}{26}} (-1+O(q))
   +2 \zeta^{4 \nu} q^{\frac{47}{26}} (-1+O(q))+\\
  &+\zeta^{2 \nu} q^{\frac{17}{26}} (6+O(q))+\zeta^{5 \nu}
   q^{\frac{23}{26}} (8+O(q))
   +\zeta^{6 \nu} q^{\frac{25}{26}} (-1+O(q))+\\
  &+\zeta^{11 \nu} q^{\frac{9}{26}} (-3+O(q))+\zeta^{8 \nu}
   q^{\frac{29}{26}} (2+O(q))
   +\zeta^{7 \nu} q^{\frac{1}{26}} (1+O(q))]^2\\
  &\times [-13 q^{\frac{7}{4}} (1+O(q))+\\
  &+\zeta^{\nu} q^{\frac{43}{52}} (13+O(q))+\zeta^{2 \nu}
   q^{\frac{47}{52}} (-22+O(q))
   +\zeta^{3 \nu} q^{\frac{51}{52}} (-21+O(q))+\\
  &+\zeta^{4 \nu} q^{\frac{3}{52}} (-1+O(q))+\zeta^{5 \nu}
   q^{\frac{7}{52}} (2+O(q))
   +\zeta^{6 \nu} q^{\frac{11}{52}} (2+O(q))+\\
  &+\zeta^{7 \nu} q^{\frac{15}{52}} (-2+O(q))+\zeta^{8 \nu}
   q^{\frac{19}{52}} (-8+O(q))
   +\zeta^{9 \nu} q^{\frac{23}{52}} (6+O(q))+\\
  &+\zeta^{10 \nu} q^{\frac{27}{52}} (1+O(q))+\zeta^{11 \nu}
   q^{\frac{31}{52}} (-8+O(q))
   +\zeta^{12 \nu} q^{\frac{35}{52}} (17+O(q))]^2
\endaligned$$
We will calculate the $q^{\frac{1}{2}}$-term which is the lowest
degree, there are five such terms:

\noindent (1)
$$\aligned
 &(\zeta^{7 \nu} q^{\frac{1}{26}} \cdot 1)^2 \times
  [(\zeta^{6 \nu} q^{\frac{11}{52}} \cdot 2)^2+
  2 \cdot \zeta^{5 \nu} q^{\frac{7}{52}} \cdot 2 \cdot
  \zeta^{7 \nu} q^{\frac{15}{52}} \cdot (-2)+\\
 &+2 \cdot \zeta^{4 \nu} q^{\frac{3}{52}} \cdot (-1) \cdot
  \zeta^{8 \nu} q^{\frac{19}{52}} \cdot (-8)]\\
=&12 q^{\frac{1}{2}}.
\endaligned$$
\noindent (2)
$$2 \cdot 2 \zeta^{9 \nu} q^{\frac{5}{26}} \cdot 1 \cdot \zeta^{7 \nu}
  q^{\frac{1}{26}} \cdot 1
  \times [(\zeta^{5 \nu} q^{\frac{7}{52}} \cdot 2)^2+
  2 \zeta^{4 \nu} q^{\frac{3}{52}} \cdot (-1) \cdot \zeta^{6 \nu}
  q^{\frac{11}{52}} \cdot 2]=0.$$
\noindent (3)
$$2 \cdot 2 \zeta^{10 \nu} q^{\frac{7}{26}} \cdot (-1) \cdot
  \zeta^{7 \nu} q^{\frac{1}{26}} \cdot 1 \times 2 \zeta^{4 \nu}
  q^{\frac{3}{52}} \cdot (-1) \cdot \zeta^{5 \nu} q^{\frac{7}{52}} \cdot 2
 =16 q^{\frac{1}{2}}.$$
\noindent (4)
$$2 \cdot \zeta^{11 \nu} q^{\frac{9}{26}} \cdot (-3) \cdot \zeta^{7 \nu}
  q^{\frac{1}{26}} \cdot 1 \times [\zeta^{4 \nu} q^{\frac{3}{52}} \cdot (-1)]^2
 =-6 q^{\frac{1}{2}}.$$
\noindent (5)
$$(2 \zeta^{9 \nu} q^{\frac{5}{26}} \cdot 1)^2 \times
  (\zeta^{4 \nu} q^{\frac{3}{52}} \cdot (-1))^2=4 q^{\frac{1}{2}}.$$
Hence, for $\Phi_{2, 2}(u_1(z), \ldots, u_6(z))$ which is a modular
form for $\Gamma(1)$ with weight $100$, the lowest degree term is given by
$$(12+0+16-6+4) q^{\frac{1}{2}} \cdot q^{\frac{9}{24} \cdot 20}=26 q^8.$$
Thus,
$$\Phi_{2, 2}(u_1(z), \ldots, u_6(z))=q^8 (13 \cdot 26+O(q)).$$
The leading term of $\Phi_{2, 2}(u_1(z), \ldots, u_6(z))$ together
with its weight $100$ suffice to identify this modular form with
$\Phi_{2, 2}(u_1(z), \ldots, u_6(z))=13 \cdot 26 \Delta(z)^8 E_4(z)$.
Consequently,
$$\Phi_{2, 2}(x_1(z), \ldots, x_6(z))=13 \cdot 26 \Delta(z)^8
  E_4(z)/\eta(z)^{160}=13 \cdot 26 \eta(z)^8 \Delta(z) E_4(z).$$
The calculation of $\Phi_{3, 3}$ and $\Phi_{6, 1}$ is similar as
that of $\Phi_{2, 2}$. We have
$$\Phi_{3, 3}(x_1(z), \ldots, x_6(z))=-13 \cdot 27 \Delta(z)^2 E_6(z).$$
$$\Phi_{6, 1}(x_1(z), \ldots, x_6(z))=-13 \cdot 285 \Delta(z)^2 E_6(z).$$
After normalization, this completes the proof of Theorem 3.3.

\noindent
$\qquad \qquad \qquad \qquad \qquad \qquad \qquad \qquad \qquad
 \qquad \qquad \qquad \qquad \qquad \qquad \qquad \boxed{}$

  Now we give the normalization for the following three families of
invariant homogeneous polynomials of degrees $d=18$, $22$ and $34$.
For $d=18$, there are two such invariant homogeneous polynomials
$\Phi_{3, 1}$ and $\Phi_{0, 3}$:
$$\left\{\aligned
  \Phi_{3, 1} &:=\frac{1}{13 \cdot 2} \left(\sum_{\nu=0}^{12}
  w_{\nu}^3 \delta_{\nu} +w_{\infty}^3 \delta_{\infty}\right),\\
  \Phi_{0, 3} &:=\frac{1}{13 \cdot 6} \left(\sum_{\nu=0}^{12}
                 \delta_{\nu}^3+\delta_{\infty}^3\right),
\endaligned\right.\eqno{(3.22)}$$
For $d=22$, there are two such invariant homogeneous polynomials
$\Phi_{1, 3}$ and $\Phi_{4, 1}$. For $d=34$, there are three such
invariant homogeneous polynomials $\Phi_{1, 5}$, $\Phi_{4, 3}$ and
$\Phi_{7, 1}$.

\textbf{Theorem 3.4.} {\it The $G$-invariant homogeneous polynomials
$\Phi_{m, n}$ of degrees $d=18$, $22$ and $34$ in $x_1(z)$, $\ldots$,
$x_6(z)$ can be identified with modular forms as follows$:$
$$\left\{\aligned
  \Phi_{3, 1}(x_1(z), \ldots, x_6(z)) &=\Delta(z) E_6(z),\\
  \Phi_{0, 3}(x_1(z), \ldots, x_6(z)) &=\Delta(z) E_6(z),\\
  \Phi_{1, 3}(x_1(z), \ldots, x_6(z)) &=0,\\
  \Phi_{4, 1}(x_1(z), \ldots, x_6(z)) &=0,\\
  \Phi_{1, 5}(x_1(z), \ldots, x_6(z)) &=0,\\
  \Phi_{4, 3}(x_1(z), \ldots, x_6(z)) &=0,\\
  \Phi_{7, 1}(x_1(z), \ldots, x_6(z)) &=0.\\
\endaligned\right.\eqno{(3.23)}$$}

{\it Proof}. The proof is similar as that of Theorem 3.3.

\noindent
$\qquad \qquad \qquad \qquad \qquad \qquad \qquad \qquad \qquad
 \qquad \qquad \qquad \qquad \qquad \qquad \qquad \boxed{}$

  As a consequence of Theorem 3.3 and Theorem 3.4, we find an image
of the modular curve $X=X(13)$ in $\mathbb{CP}^5$ given by the complete
intersection curve $Y$.

\textbf{Theorem 3.5.} {\it There is a morphism
$$\Phi: X \to Y \subset \mathbb{CP}^5$$
with $\Phi(z)=(x_1(z), \ldots, x_6(z))$, where $Y$ is a complete
intersection curve corresponding to the ideal
$$I=I(Y)=(\Phi_{4}, \Phi_{8}, \Phi_{10}, \Phi_{14}).\eqno{(3.24)}$$
Moreover, there is a morphism $\Psi: X \to Z \subset \mathbb{CP}^5$,
where $Z$ is the scheme corresponding to the ideal
$$J=J(Z)=(\Phi_4, \Phi_8, \Phi_{10}, \Phi_{14}, \Phi_{4, 0}, \Phi_{1, 2},
         \Phi_{1, 3}, \Phi_{4, 1}, \Phi_{1, 5}, \Phi_{4, 3}, \Phi_{7, 1}).
         \eqno{(3.25)}$$}

{\it Proof}. Theorem 3.3 and Theorem 3.4 imply that
$$\left\{\aligned
  \Phi_{4}(x_1(z), \ldots, x_6(z)) &=0,\\
  \Phi_{8}(x_1(z), \ldots, x_6(z)) &=0,\\
  \Phi_{10}(x_1(z), \ldots, x_6(z)) &=0,\\
  \Phi_{14}(x_1(z), \ldots, x_6(z)) &=0,\\
  \Phi_{4, 0}(x_1(z), \ldots, x_6(z)) &=0,\\
  \Phi_{1, 2}(x_1(z), \ldots, x_6(z)) &=0,\\
  \Phi_{1, 3}(x_1(z), \ldots, x_6(z)) &=0,\\
  \Phi_{4, 1}(x_1(z), \ldots, x_6(z)) &=0,\\
  \Phi_{1, 5}(x_1(z), \ldots, x_6(z)) &=0,\\
  \Phi_{4, 3}(x_1(z), \ldots, x_6(z)) &=0,\\
  \Phi_{7, 1}(x_1(z), \ldots, x_6(z)) &=0.
\endaligned\right.$$
This complete the proof of Theorem 3.5.

\noindent
$\qquad \qquad \qquad \qquad \qquad \qquad \qquad \qquad \qquad
 \qquad \qquad \qquad \qquad \qquad \qquad \qquad \boxed{}$

\begin{center}
{\large\bf 4. A different construction: $E_8$-singularity from
              $C_Y/\text{SL}(2, 13)$ over $X(13)$ and a variation
              of the $E_8$-singularity structure over $X(13)$}
\end{center}

  In this section, we will give a different construction of the
$E_8$-singularity from a quotient $C_Y/G$ over the modular curve $X$.
The significance of the complete intersection curve $Y$ with
multi-degree $(4, 8, 10, 14)$ is that the finite group $G$ acts
linearly on $\mathbb{C}^6$ and on $\mathbb{CP}^5$ leaving invariant
$Y \subset \mathbb{CP}^5$ and the cone $C_Y \subset \mathbb{C}^6$.

\textbf{Theorem 4.1.} (A variation of the $E_8$-singularity structure
over the modular curve $X$: algebraic version) {\it The equation of
$E_8$-singularity $$\Phi_{20}^3-\Phi_{30}^2-1728 \Phi_{12}^5=0$$
possesses an infinitely many kinds of distinct modular parametrizations
$($with the cardinality of the continuum in ZFC set theory$)$
$$(\Phi_{12}, \Phi_{20}, \Phi_{30})=(\Phi_{12}^{\lambda},
   \Phi_{20}^{\mu}, \Phi_{30}^{\gamma})\eqno{(4.1)}$$
over the modular curve $X$ as follows$:$
$$\left\{\aligned
  \Phi_{12}^{\lambda} &=\lambda \Phi_{3, 0}+(1-\lambda) \Phi_{0, 2}
              \quad \text{mod $\mathfrak{a}_1$},\\
  \Phi_{20}^{\mu} &=\mu \Phi_{5, 0}+(1-\mu) \Phi_{2, 2}
              \quad \text{mod $\mathfrak{a}_2$},\\
  \Phi_{30}^{\gamma} &=\gamma_1 \Phi_{0, 5}+\gamma_2 \Phi_{3, 3}+
              (1-\gamma_1-\gamma_2) \Phi_{6, 1} \quad \text{mod $\mathfrak{a}_3$},
\endaligned\right.\eqno{(4.2)}$$
where $\Phi_{12}$, $\Phi_{20}$ and $\Phi_{30}$ are invariant homogeneous
polynomials of degree $12$, $20$ and $30$, respectively. The ideals are
given by
$$\left\{\aligned
  \mathfrak{a}_1 &=(\Phi_4, \Phi_8),\\
  \mathfrak{a}_2 &=(\Phi_4, \Phi_8, \Phi_{10}, \Phi_{4, 0}, \Phi_{1, 2}),\\
  \mathfrak{a}_3 &=(\Phi_4, \Phi_8, \Phi_{10}, \Phi_{3, 0}, \Phi_{0, 2},
                    \Phi_{14}, \Phi_{4, 0}, \Phi_{1, 2}, \Phi_{1, 3},
                    \Phi_{4, 1}),
\endaligned\right.\eqno{(4.3)}$$
and the parameter space $\{ (\lambda, \mu, \gamma) \} \cong \mathbb{C}^4$.
They form a variation of the $E_8$-singularity structure over the modular
curve $X$.}

{\it Proof}. By Theorem 3.5, the ideals $\mathfrak{a}_1$ and $\mathfrak{a}_2$
are zero ideal over the modular curve $X$, and $\mathfrak{a}_3=(\Delta(z))$
over the modular curve $X$. On the other hand, by Theorem 3.3, for degree
$d=12$, the invariant homogeneous polynomials $\Phi_{3, 0}$ and $\Phi_{0, 2}$
form a two-dimensional complex vector space, and
$\Phi_{3, 0}=\Phi_{0, 2}=\Delta(z)$ over the modular curve $X$ (after
normalization). Hence, $\Phi_{12}^{\lambda}=\Delta(z)$ over the modular curve
$X$. For degree $d=20$, the invariant homogeneous polynomials $\Phi_{5, 0}$
and $\Phi_{2, 2}$ form a two-dimensional complex vector space, and
$\Phi_{5, 0}=\Phi_{2, 2}=\eta(z)^8 \Delta(z) E_4(z)$ over the modular curve $X$
(after normalization). Hence, $\Phi_{20}^{\mu}=\eta(z)^8 \Delta(z) E_4(z)$ over
the modular curve $X$. Finally, for degree $d=30$, the invariant homogeneous
polynomials $\Phi_{0, 5}$, $\Phi_{3, 3}$ and $\Phi_{6, 1}$ form a
three-dimensional complex vector space, and
$\Phi_{0, 5}=\Phi_{3, 3}=\Phi_{6, 1}=\Delta(z)^2 E_6(z)$ over the modular curve
$X$ (after normalization). Hence, $\Phi_{30}^{\gamma}=\Delta(z)^2 E_6(z)$ over
the modular curve $X$. This shows that
$$\Phi_{20}^3-\Phi_{30}^2-1728 \Phi_{12}^5=0$$
for $(\Phi_{12}, \Phi_{20}, \Phi_{30})=(\Phi_{12}^{\lambda}, \Phi_{20}^{\mu},
\Phi_{30}^{\gamma})$ over $X$.

\noindent
$\qquad \qquad \qquad \qquad \qquad \qquad \qquad \qquad \qquad
 \qquad \qquad \qquad \qquad \qquad \qquad \qquad \boxed{}$

  In order to define a $G$-invariant map from $C_Y$ to the $E_8$-singularity,
the morphism $\Phi: X \to Y \subset \mathbb{CP}^5$ given by
$\Phi(z)=(x_1(z), \ldots, x_6(z))$ should be revised as follows: near
a cusp, for this to extend at the cusp, $x_i(z)$ has to be replaced by
$q^{-\frac{2}{3 \cdot 13}} x_i(z)$ for $1 \leq i \leq 6$.

\textbf{Theorem 4.2.} (A variation of the $E_8$-singularity structure over
the modular curve $X$: geometric version) {\it There is a morphism from the
cone $C_Y$ over $Y$ to the $E_8$-singularity:
$$f: C_Y/G \rightarrow \text{Spec} \left(\mathbb{C}[\Phi_{12}, \Phi_{20},
  \Phi_{30}]/(\Phi_{20}^3-\Phi_{30}^2-1728 \Phi_{12}^5)\right) \eqno{(4.4)}$$
over the modular curve $X$. In particular, there are infinitely many such
triples $(\Phi_{12}, \Phi_{20}, \Phi_{30})$ $=$ $(\Phi_{12}^{\lambda}, \Phi_{20}^{\mu},
\Phi_{30}^{\gamma})$ whose parameter space $\{ (\lambda, \mu, \gamma) \}$
$\cong \mathbb{C}^4$. They form a variation of the $E_8$-singularity structure
over the modular curve $X$.}

{\it Proof}. By Theorem 3.3, Theorem 3.5 and Theorem 4.1, we have the following
ring homomorphism
$$\mathbb{C}[\Phi_{12}, \Phi_{20}, \Phi_{30}]/(\Phi_{20}^3-\Phi_{30}^2-
  1728 \Phi_{12}^5) \rightarrow \left[\mathbb{C}[z_1, z_2, z_3, z_4, z_5,
  z_6]/I\right]^{\text{SL}(2, 13)}$$
over the modular curve $X$, which induces a morphism of schemes
$$C_Y/G \rightarrow \text{Spec} \left(\mathbb{C}[\Phi_{12}, \Phi_{20},
  \Phi_{30}]/(\Phi_{20}^3-\Phi_{30}^2-1728 \Phi_{12}^5)\right)$$
over $X$. This completes the proof of Theorem 4.2.

\noindent
$\qquad \qquad \qquad \qquad \qquad \qquad \qquad \qquad \qquad
 \qquad \qquad \qquad \qquad \qquad \qquad \qquad \boxed{}$

\textbf{Theorem 4.3.} (A variation of the structure of decomposition
formulas of the elliptic modular functions $j$ over the modular curve
$X$) {\it There are infinitely many kinds of distinct decomposition
formulas of the elliptic modular function $j$ in terms of the invariants
$\Phi_{12}$, $\Phi_{20}$ and $\Phi_{30}$ over the modular curve $X$:
$$j(z): j(z)-1728: 1=\Phi_{20}^3: \Phi_{30}^2: \Phi_{12}^5,
  \eqno{(4.5)}$$
where $(\Phi_{12}, \Phi_{20}, \Phi_{30})=(\Phi_{12}^{\lambda},
\Phi_{20}^{\mu}, \Phi_{30}^{\gamma})$ are given by (4.2). They form a
variation of the structure of decomposition formulas of the elliptic
modular functions $j$ over the modular curve $X$.}

{\it Proof}. By Theorem 3.3, we have the following relations over the
modular curve $X$:
$$j(z)=\frac{E_4(z)^3}{\Delta(z)}=\frac{\Phi_{20}^3}{\Phi_{12}^5}, \quad
  j(z)-1728=\frac{E_6(z)^2}{\Delta(z)}=\frac{\Phi_{30}^2}{\Phi_{12}^5}.
  \eqno{(4.6)}$$
Hence, we obtain an infinitely many kinds of distinct decomposition formulas
of the elliptic modular function $j$ in terms of the invariants $\Phi_{12}$,
$\Phi_{20}$ and $\Phi_{30}$ over the modular curve $X$:
$$j(z): j(z)-1728: 1=\Phi_{20}^3: \Phi_{30}^2: \Phi_{12}^5,
  \eqno{(4.7)}$$
where $(\Phi_{12}, \Phi_{20}, \Phi_{30})=(\Phi_{12}^{\lambda}, \Phi_{20}^{\mu},
\Phi_{30}^{\gamma})$ are given by (4.2). This completes the proof of Theorem
4.3.

\noindent
$\qquad \qquad \qquad \qquad \qquad \qquad \qquad \qquad \qquad
 \qquad \qquad \qquad \qquad \qquad \qquad \qquad \boxed{}$

  Theorem 4.3 shows that there are infinitely many kinds of distinct
decomposition formulas of the elliptic modular function $j$ in terms of the
invariant polynomials of the same degrees $12$, $20$ and $30$ over the modular
curves, one and only one is given by $\text{SL}(2, 5)$ corresponding to the
modular curve $X(5)$, the other infinitely many kinds of distinct decomposition
formulas (which form a variation of the structure of decomposition formulas)
are given by $\text{SL}(2, 13)$ corresponding to the modular curve $X(13)$.

  In the end, let us recall some facts about exotic spheres (see \cite{Hi}).
A $k$-dimensional compact oriented differentiable manifold is called a
$k$-sphere if it is homeomorphic to the $k$-dimensional standard sphere.
A $k$-sphere not diffeomorphic to the standard $k$-sphere is said to be
exotic. The first exotic sphere was discovered by Milnor in 1956 (see
\cite{Mi}). Two $k$-spheres are called equivalent if there exists an
orientation preserving diffeomorphism between them. The equivalence classes
of $k$-spheres constitute for $k \geq 5$ a finite abelian group $\Theta_k$
under the connected sum operation. $\Theta_k$ contains the subgroup $b P_{k+1}$
of those $k$-spheres which bound a parallelizable manifold. $b P_{4m}$ ($m \geq 2$)
is cyclic of order $2^{2m-2}(2^{2m-1}-1)$ numerator $(4 B_m/m)$, where $B_m$
is the $m$-th Bernoulli number. Let $g_m$ be the Milnor generator of $b P_{4m}$.
If a $(4m-1)$-sphere $\Sigma$ bounds a parallelizable manifold $B$ of dimension
$4m$, then the signature $\tau(B)$ of the intersection form of $B$ is divisible
by $8$ and $\Sigma=\frac{\tau(B)}{8} g_m$. For $m=2$ we have
$b P_8=\Theta_7=\mathbb{Z} /28 \mathbb{Z}$. All these results are due to
Milnor-Kervaire (see \cite{KM}). In particular,
$$\sum_{i=0}^{2m} z_i \overline{z_i}=1, \quad
  z_0^3+z_1^{6k-1}+z_2^2+\cdots+z_{2m}^2=0$$
is a $(4m-1)$-sphere embedded in $S^{4m+1} \subset \mathbb{C}^{2n+1}$ which
represents the element $(-1)^m k \cdot g_m \in b P_{4m}$. For $m=2$ and
$k=1, 2, \cdots, 28$ we get the $28$ classes of $7$-spheres. Theorem 4.1 and
Theorem 4.2 show that the higher dimensional liftings of infinitely many kinds
of distinct constructions of the $E_8$-singularity: $\mathbb{C}^2/\text{SL}(2, 5)$
and a variation of the $E_8$-singularity structure
$$C_Y/\text{SL}(2, 13) \rightarrow \text{Spec} \left(\mathbb{C}[\Phi_{12},
  \Phi_{20}, \Phi_{30}]/(\Phi_{20}^3-\Phi_{30}^2-1728 \Phi_{12}^5)\right)$$
over the modular curve $X$ give the same Milnor's standard generator
of $\Theta_7$. This is the version of differential topology.

\begin{center}
{\large\bf 5. $Q_{18}$ and $E_{20}$-singularities from $C_Y/\text{SL}(2, 13)$
              over $X(13)$ and variations of $Q_{18}$ and $E_{20}$-singularity
              structures over $X(13)$}
\end{center}

  In this section, we will extend our work on a variation of the
$E_8$-singularity structure over the modular curve $X$ to the cases
of $Q_{18}$ and $E_{20}$-singularities and obtain variations of 
$Q_{18}$ and $E_{20}$-singularity structures over the modular curve 
$X$.

  Now we give the normalization for the following three families of
invariant homogeneous polynomials of degrees $d=32$, $42$ and $44$.
For $d=32$, there are three such invariant homogeneous polynomials
$\Phi_{8, 0}$, $\Phi_{5, 2}$ and $\Phi_{2, 4}$:
$$\left\{\aligned
  \Phi_{8, 0} &:=-\frac{1}{13 \cdot 1840} \left(\sum_{\nu=0}^{12} w_{\nu}^8
                 +w_{\infty}^8\right),\\
  \Phi_{5, 2} &:=-\frac{1}{13 \cdot 2064} \left(\sum_{\nu=0}^{12} w_{\nu}^5
                 \delta_{\nu}^2+w_{\infty}^5 \delta_{\infty}^2\right),\\
  \Phi_{2, 4} &:=-\frac{1}{13 \cdot 680} \left(\sum_{\nu=0}^{12} w_{\nu}^2
                 \delta_{\nu}^4+w_{\infty}^2 \delta_{\infty}^4\right).
\endaligned\right.\eqno{(5.1)}$$
For $d=42$, there are four such invariant homogeneous polynomials
$\Phi_{0, 7}$, $\Phi_{3, 5}$, $\Phi_{6, 3}$ and $\Phi_{9, 1}$:
$$\left\{\aligned
  \Phi_{0, 7} &:=\frac{1}{13 \cdot 226842} \left(\sum_{\nu=0}^{12}
                 \delta_{\nu}^7+\delta_{\infty}^7\right),\\
  \Phi_{3, 5} &:=\frac{1}{13 \cdot 634} \left(\sum_{\nu=0}^{12} w_{\nu}^3
                 \delta_{\nu}^5+w_{\infty}^3 \delta_{\infty}^5\right),\\
  \Phi_{6, 3} &:=\frac{1}{13 \cdot 10656} \left(\sum_{\nu=0}^{12} w_{\nu}^6
                 \delta_{\nu}^3+w_{\infty}^6 \delta_{\infty}^3\right),\\
  \Phi_{9, 1} &:=\frac{1}{13 \cdot 39134} \left(\sum_{\nu=0}^{12} w_{\nu}^9
                 \delta_{\nu}+w_{\infty}^9 \delta_{\infty}\right).
\endaligned\right.\eqno{(5.2)}$$
For $d=44$, there are four such invariant homogeneous polynomials
$\Phi_{11, 0}$, $\Phi_{8, 2}$, $\Phi_{5, 4}$ and $\Phi_{2, 6}$. Let
$$\Phi_{11, 0}:=\frac{1}{13 \cdot 146905} \left(\sum_{\nu=0}^{12} w_{\nu}^{11}
                +w_{\infty}^{11}\right).\eqno{(5.3)}$$

\textbf{Theorem 5.1.} {\it The $G$-invariant homogeneous polynomials
$\Phi_{m, n}$ of degrees $d=32$, $42$ and $44$ in $x_1(z)$, $\ldots$,
$x_6(z)$ can be identified with modular forms as follows$:$
$$\left\{\aligned
  \Phi_{8, 0}(x_1(z), \ldots, x_6(z)) &=\eta(z)^8 \Delta(z)^2 E_4(z),\\
  \Phi_{5, 2}(x_1(z), \ldots, x_6(z)) &=\eta(z)^8 \Delta(z)^2 E_4(z),\\
  \Phi_{2, 4}(x_1(z), \ldots, x_6(z)) &=\eta(z)^8 \Delta(z)^2 E_4(z),\\
  \Phi_{0, 7}(x_1(z), \ldots, x_6(z)) &=\Delta(z)^3 E_6(z),\\
  \Phi_{3, 5}(x_1(z), \ldots, x_6(z)) &=\Delta(z)^3 E_6(z),\\
  \Phi_{6, 3}(x_1(z), \ldots, x_6(z)) &=\Delta(z)^3 E_6(z),\\
  \Phi_{9, 1}(x_1(z), \ldots, x_6(z)) &=\Delta(z)^3 E_6(z),\\
  \Phi_{11, 0}(x_1(z), \ldots, x_6(z)) &=\eta(z)^8 \Delta(z)^3 E_4(z),\\
  \Phi_{8, 2}(x_1(z), \ldots, x_6(z)) &\in \eta(z)^8 \Delta(z)^2
  (\mathbb{C} E_4(z)^4 \oplus \mathbb{C} E_4(z) E_6(z)^2),\\
  \Phi_{5, 4}(x_1(z), \ldots, x_6(z)) &\in \eta(z)^8 \Delta(z)^2
  (\mathbb{C} E_4(z)^4 \oplus \mathbb{C} E_4(z) E_6(z)^2),\\
  \Phi_{2, 6}(x_1(z), \ldots, x_6(z)) &\in \eta(z)^8 \Delta(z)^2
  (\mathbb{C} E_4(z)^4 \oplus \mathbb{C} E_4(z) E_6(z)^2).
\endaligned\right.\eqno{(5.4)}$$}

{\it Proof}. The proof is similar as that of Theorem 3.3. For the
calculation of $\Phi_{8, 0}$, $\Phi_{0, 7}$ and $\Phi_{11, 0}$,
see also \cite{Y3}, section 4.

\noindent
$\qquad \qquad \qquad \qquad \qquad \qquad \qquad \qquad \qquad
 \qquad \qquad \qquad \qquad \qquad \qquad \qquad \boxed{}$

\textbf{Theorem 5.2.} (Variations of $Q_{18}$ and $E_{20}$-singularity
structures over the modular curve $X$: algebraic version) {\it The equations
of $Q_{18}$ and $E_{20}$-singularities
$$\Phi_{32}^3-\Phi_{12} \Phi_{42}^2-1728 \Phi_{12}^8=0, \quad
  \Phi_{44}^3-\Phi_{12}^4 \Phi_{42}^2-1728 \Phi_{12}^{11}=0$$
possess an infinitely many kinds of distinct modular parametrizations
$($with the cardinality of the continuum in ZFC set theory$)$
$$(\Phi_{12}, \Phi_{32}, \Phi_{42}, \Phi_{44})=(\Phi_{12}^{\lambda},
   \Phi_{32}^{\mu}, \Phi_{42}^{\gamma}, \Phi_{44})\eqno{(5.5)}$$
over the modular curve $X$ as follows$:$
$$\left\{\aligned
  \Phi_{12}^{\lambda} &=\lambda \Phi_{3, 0}+(1-\lambda) \Phi_{0, 2}
              \quad \text{mod $\mathfrak{a}_1$},\\
  \Phi_{32}^{\mu} &=\mu_1 \Phi_{8, 0}+\mu_2 \Phi_{5, 2}+(1-\mu_1-\mu_2)
              \Phi_{2, 4} \quad \text{mod $\mathfrak{a}_3$},\\
  \Phi_{42}^{\gamma} &=\gamma_1 \Phi_{0, 7}+\gamma_2 \Phi_{3, 5}+\gamma_3
              \Phi_{6, 3}+(1-\gamma_1-\gamma_2-\gamma_3) \Phi_{9, 1}
              \quad \text{mod $\mathfrak{a}_4$},\\
  \Phi_{44} &=\Phi_{11, 0} \quad \text{mod $\mathfrak{a}_4$},
\endaligned\right.\eqno{(5.6)}$$
where $\Phi_{12}$, $\Phi_{32}$, $\Phi_{42}$ and $\Phi_{44}$ are invariant
homogeneous polynomials of degree $12$, $32$, $42$ and $44$, respectively.
The ideals are given by
$$\left\{\aligned
  \mathfrak{a}_1 &=(\Phi_4, \Phi_8),\\
  \mathfrak{a}_3 &=(\Phi_4, \Phi_8, \Phi_{10}, \Phi_{3, 0}, \Phi_{0, 2},
  \Phi_{14}, \Phi_{4, 0}, \Phi_{1, 2}, \Phi_{1, 3}, \Phi_{4, 1}),\\
  \mathfrak{a}_4 &=(\Phi_4, \Phi_8, \Phi_{10}, \Phi_{3, 0}, \Phi_{0, 2},
  \Phi_{14}, \Phi_{4, 0}, \Phi_{1, 2}, \Phi_{1, 3}, \Phi_{4, 1}, \Phi_{1, 5},
  \Phi_{4, 3}, \Phi_{7, 1}),
\endaligned\right.\eqno{(5.7)}$$
and the parameter space $\{ (\lambda, \mu, \gamma) \} \cong \mathbb{C}^6$.
They form variations of $Q_{18}$ and $E_{20}$-singularity structures over
the modular curve $X$.}

{\it Proof}. The proof is similar as that of Theorem 4.1.

\noindent
$\qquad \qquad \qquad \qquad \qquad \qquad \qquad \qquad \qquad
 \qquad \qquad \qquad \qquad \qquad \qquad \qquad \boxed{}$

  In order to define $G$-invariant maps from $C_Y$ to $Q_{18}$ and
$E_{20}$-singularities, the morphism $\Phi: X \to Y \subset \mathbb{CP}^5$
given by $\Phi(z)=(x_1(z), \ldots, x_6(z))$ should be revised as follows:
near a cusp, for this to extend at the cusp, $x_i(z)$ has to be replaced
by $q^{-\frac{2}{3 \cdot 13}} x_i(z)$ for $1 \leq i \leq 6$.

\textbf{Theorem 5.3.} (Variations of $Q_{18}$ and $E_{20}$-singularity
structures over the modular curve $X$: geometric version) {\it There are two
morphisms from the cone $C_Y$ over $Y$ to the $Q_{18}$ and $E_{20}$-singularities:
$$f_1: C_Y/G \rightarrow \text{Spec} \left(\mathbb{C}[\Phi_{12}, \Phi_{32},
           \Phi_{42}]/(\Phi_{32}^3-\Phi_{12} \Phi_{42}^2-1728 \Phi_{12}^8)\right)
           \eqno{(5.8)}$$
and
$$f_2: C_Y/G \rightarrow \text{Spec} \left(\mathbb{C}[\Phi_{12}, \Phi_{42},
           \Phi_{44}]/(\Phi_{44}^3-\Phi_{12}^4 \Phi_{42}^2-1728 \Phi_{12}^{11})\right)
           \eqno{(5.9)}$$
over the modular curve $X$. In particular, there are infinitely many such
triples $(\Phi_{12}, \Phi_{32}, \Phi_{42})$ $=$ $(\Phi_{12}^{\lambda}, \Phi_{32}^{\mu},
\Phi_{42}^{\gamma})$ whose parameter space $\{ (\lambda, \mu, \gamma) \}$
$\cong \mathbb{C}^6$. They form variations of $Q_{18}$ and $E_{20}$-singularity
structures over the modular curve $X$.}

{\it Proof}. The proof is similar as that of Theorem 4.2.

\noindent
$\qquad \qquad \qquad \qquad \qquad \qquad \qquad \qquad \qquad
 \qquad \qquad \qquad \qquad \qquad \qquad \qquad \boxed{}$

  In fact, Theorem 5.1 and Theorem 5.2 show that the higher dimensional
liftings of infinitely many kinds of distinct constructions of the
$E_{20}$-singularity, i.e., a variation of the $E_{20}$-singularity structure
$$C_Y/\text{SL}(2, 13) \rightarrow \text{Spec} \left(\mathbb{C}
  [\Phi_{12}, \Phi_{42}, \Phi_{44}]/(\Phi_{44}^3-\Phi_{12}^4
  \Phi_{42}^2-1728 \Phi_{12}^{11})\right)$$
over the modular curve $X$ give the square of the Milnor's standard
generator of $\Theta_7$:
$$z_1^{11}+z_2^3+z_3^2+z_4^2+z_5^2=0, \quad \sum_{i=1}^{5} z_i
  \overline{z_i}=1, \quad z_i \in \mathbb{C} \quad (1 \leq i \leq 5).$$

\vskip 2.0 cm

\noindent{Department of Mathematics, Peking University}

\noindent{Beijing 100871, P. R. China}

\noindent{\it E-mail address}: yanglei$@$math.pku.edu.cn


\vskip 1.5 cm

\begin{thebibliography}{99}

\bibitem{Ar1} V. I. Arnol'd, Critical points of smooth functions,
              in: {\it Proceedings of the International Congress
              of Mathematicians $($Vancouver, B. C., 1974$)$},
              Vol. 1, 19-39, Canad. Math. Congress, Montreal,
              Que., 1975.

\bibitem{Ar2} V. I. Arnol'd, S. M. Gusein-Zade and  A. N. Varchenko,
             {\it Singularities of Differentiable Maps}, Vol. I.
             {\it The Classification of Critical Points, Caustics
             and Wave Fronts}, Translated from the Russian by Ian
             Porteous and Mark Reynolds, Monographs in Mathematics,
             {\bf 82}, Birkh\"{a}user, 1985.

\bibitem{Ar3} V. I. Arnol'd, {\it Catastrophe Theory}, translated
             from the Russian by G. S. Wassermann, based on a
             translation by R. K. Thomas, Third edition, Springer-Verlag,
             Berlin, 1992.

\bibitem{Br1} E. Brieskorn, Beispiele zur Differentialtopologie von
             Singularit\"{a}ten, Invent. Math. {\bf 2} (1966), 1-14.

\bibitem{Br2} E. Brieskorn, Singular elements of semi-simple algebraic
              groups, in: {\it Actes du Congr\`{e}s International des
              Math\'{e}maticiens $($Nice, 1970$)$}, Tome 2, 279-284,
              Gauthier-Villars, Paris, 1971.

\bibitem{Br3} E. Brieskorn, Singularities in the work of Friedrich
              Hirzebruch, in: {\it Surveys in differential geometry},
              17-60, Surv. Differ. Geom., {\bf 7}, Int. Press,
              Somerville, MA, 2000.

\bibitem{BrPR} E. Brieskorn, A. Pratoussevitch and F. Rothenh\"{a}usler,
               The combinatorial geometry of singularities and Arnold's
               series $E$, $Z$, $Q$, Dedicated to Vladimir I. Arnold on
               the occasion of his 65th birthday, Mosc. Math. J. {\bf 3}
               (2003), 273-333.

\bibitem{Do} I. V. Dolgachev, Quotient-conical singularities on complex
             surfaces, Funktsional. Anal. i Prilozhen. {\bf 8} (1974),
             75-76 (Russian), English translation in: Funct. Anal. Appl.
             {\bf 8} (1974), 160-161.

\bibitem{Du} W. Duke, Continued fractions and modular functions, Bull.
             Amer. Math. Soc. (N.S.) {\bf 42} (2005), 137-162.

\bibitem{FK} H. M. Farkas and I. Kra, {\it Theta Constants, Riemann
            Surfaces and the Modular Group, An Introduction with
            Applications to Uniformization Theorems, Partition
            Identities and Combinatorial Number Theory}, Graduate
            Studies in Mathematics, {\bf 37}, American Mathematical
            Society, Providence, RI, 2001.

\bibitem{Gr} G.-M. Greuel, Some aspects of Brieskorn's mathematical
             work, {\it Singularities $($Oberwolfach, 1996$)$},
             xv-xxii, Progr. Math., {\bf 162}, Birkh\"{a}user, Basel,
             1998.

\bibitem{GrP} G.-M. Greuel and W. Purkert, Leben und Werk von Egbert
              Brieskorn (1936-2013), Jahresber. Dtsch. Math.-Ver.
              {\bf 118} (2016), 143-178.

\bibitem{GrP2} G.-M. Greuel and W. Purkert, Life and work of Egbert
               Brieskorn (1936-2013), J. Singul. {\bf 18} (2018), 1-28.

\bibitem{He} C. Hermite, Sur la r\'{e}solution de l'\'{e}quation du
             cinqui\`{e}me degr\'{e}, C. R. Acad. Sci. Paris {\bf 46}
             (1858), 508-515, in: {\it \OE uvres de Charles Hermite},
             Vol. II, Gauthier-Villars, 1908.

\bibitem{Hi} F. Hirzebruch, Singularities and exotic spheres,
             S\'{e}minaire Bourbaki, 1966/67, Exp. 314, in: {\it
             Gesammelte Abhandlungen}, Bd. II, 70-80, Springer-Verlag,
             1987.

\bibitem{KM} M. Kervaire and J. Milnor, Groups of homotopy spheres: I,
             Ann. of Math. (2) {\bf 77} (1963), 504-537.

\bibitem{KS} R. C. Kirby and M. G. Scharlemann, Eight faces of the
             Poincar\'{e} homology $3$-sphere, in: {\it Geometric
             topology $($Proc. Georgia Topology Conf., Athens, Ga.,
             1977$)$}, 113-146, Academic Press, New York-London, 1979.

\bibitem{K} F. Klein, {\it Lectures on the Icosahedron and the
            Solution of Equations of the Fifth Degree}, Translated by
            G. G. Morrice, second and revised edition, Dover Publications,
            Inc., 1956.

\bibitem{K1} F. Klein, \"{U}ber die Transformation der elliptischen Functionen
            und die Aufl\"{o}sung der Gleichungen f\"{u}nften Grades, Math.
            Ann. {\bf 14} (1879), 111-172, in: {\it Gesammelte Mathematische
            Abhandlungen}, Bd. III, 13-75, Springer-Verlag, Berlin, 1923.

\bibitem{K2} F. Klein, \"{U}ber die Transformation siebenter Ordnung
             der elliptischen Functionen, Math. Ann. {\bf 14} (1879),
             428-471, in: {\it Gesammelte Mathematische Abhandlungen},
             Bd. III, 90-136, Springer-Verlag, Berlin, 1923.

\bibitem{KF1} F. Klein and R. Fricke, {\it Vorlesungen \"{u}ber die Theorie
            der Elliptischen Modulfunctionen}, Vol. I, Leipzig, 1890;

\bibitem{KF2} F. Klein and R. Fricke, {\it Vorlesungen \"{u}ber die Theorie
            der Elliptischen Modulfunctionen}, Vol. II, Leipzig, 1892.

\bibitem{Kr1} P. B. Kronheimer, The construction of ALE spaces as hyper-K\"{a}ahler
              quotients, J. Diff. Geom. {\bf 29} (1989), 665-683.

\bibitem{Kr2} P. B. Kronheimer, Instantons and the geometry of the nilpotent
              variety, J. Diff. Geom. {\bf 32} (1990), 473-490.

\bibitem{Mc} J. McKay, Graphs, singularities, and finite groups, in:
             {\it The Santa Cruz Conference on Finite Groups $($Univ.
             California, Santa Cruz, Calif., 1979$)$}, 183-186, Proc.
             Sympos. Pure Math., {\bf 37}, Amer. Math. Soc., Providence,
             R.I., 1980.

\bibitem{Mi} J. Milnor, On manifolds homeomorphic to the $7$-sphere,
              Ann. of Math. (2) {\bf 64} (1956), 399-405.

\bibitem{Sch} H. A. Schwarz, \"{U}ber diejenigen F\"{a}lle, in welchen
             die Gaussische hypergeometrische Reihe eine algebraische
             Function ihres vierten Elementes darstellt, J. Reine Angew.
             Math. {\bf 75} (1872), 292-335.

\bibitem{Sl1} P. Slodowy, {\it Simple Singularities and Simple Algebraic
              Groups}, Lecture Notes in Math. {\bf 815}, Springer, Berlin,
              1980.

\bibitem{Sl2} P. Slodowy, Platonic solids, Kleinian singularities, and Lie
              groups, in: {\it Algebraic geometry $($Ann Arbor, Mich., 1981$)$},
              102-138, Lecture Notes in Math., {\bf 1008}, Springer, Berlin,
              1983.

\bibitem{Y1} L. Yang, Exotic arithmetic structure on the first Hurwitz triplet,
             arXiv:1209.1783v5 [math.NT], 2013.

\bibitem{Y2} L. Yang, Dedekind $\eta$-function, Hauptmodul and invariant theory,
             arXiv:1407.3550v2 [math.NT], 2014.

\bibitem{Y3} L. Yang, Icosahedron, exceptional singularities and modular forms,
             arXiv:1511.05278 [math.NT], 2015.

\end{thebibliography}
\end{document}